\documentclass[11pt,leqno]{article}
\usepackage{graphicx, amsfonts, amsthm, amsxtra, amssymb, verbatim,latexsym}
\usepackage[mathscr]{euscript}
\usepackage[all,cmtip]{xy}

\begin{document}
\def\cy{{Calabi-Yau }}
\def\X{{\cal X}}
\def\W{{\cal W}}
\def\M{{\cal M}}
\def\DHR{{\rm DHR }}            
\def\H{{\sf H}}            
\def\Z{\mathbb{Z}}                   
\def\Q{\mathbb{Q}}                   
\def\C{\mathbb{C}}                   
\def\R{\mathbb{R}}                   
\def\N{\mathbb{N}}                   
\def\ker{{\rm ker}}              
\def\cl{{\rm cl}}                
\def\dR{{\rm dR}}                
\def\P {\mathbb{P}}                  
\def\F{{\cal F}}                     
\def\E{{\cal E}}
\def\M{{\cal M}}
\def\RR{{\cal R}}
\def\W{{\cal W}}                     
\def\Tm{{\sf T }}                    
\def\tr{{\mathsf t}}                 
\def\GM{{\sf GM}}
\def\L{{\sf L}}
\def\O{{\cal O}}
\def\b{{\sf b}}
\def\h{{\sf h}}
\newtheorem{theo}{Theorem}[section]
\newtheorem{exam}{Example}[section]
\newtheorem{coro}{Corollary}[section]
\newtheorem{defi}{Definition}[section]
\newtheorem{prob}{Problem}[section]
\newtheorem{lemm}{Lemma}[section]
\newtheorem{prop}{Proposition}[section]
\newtheorem{rem}{Remark}[section]
\newtheorem{obs}{Observation}[section]
\newtheorem{conj}{Conjecture}
\newtheorem{nota}{Notation}[section]
\newtheorem{ass}{Assumption}[section]
\newtheorem{calc}{}
\numberwithin{equation}{section}

\begin{center}
{\LARGE\bf  Darboux-Halphen-Ramanujan Vector Field on a Moduli of Calabi-Yau Manifolds}
\\
\vspace{.25in} {\large {\sc Younes Nikdelan}} \\
Instituto de Matem\'atica Pura e Aplicada (IMPA), \\
Estrada Dona Castorina, 110,\\
22460-320, Rio de Janeiro, RJ, Brazil, \\
{nikdelan@impa.br}
\end{center}

\begin{abstract}
In this paper we obtain an ordinary differential
equation $\H$ from a  Picard-Fuchs equation  associated with a nowhere vanishing
holomorphic $n$-form. We work on a moduli space $\Tm$ constructed from  a \cy $n$-fold $W$ together
with a basis of the middle complex de Rham cohomology of $W$. We verify the existence of
a unique vector field $\H$ on $\Tm$  such that its composition with the Gauss-Manin connection satisfies certain
properties. The ordinary differential equation given by
$\H$ is a generalization of differential equations introduced by Darboux, Halphen and Ramanujan.
\end{abstract}
\vspace{.25in}
{\bf Kywords} Darboux-Halphen-Ramanujan vector field. Hodge structure. Picard-Fuchs equation. Gauss-Manin connection.
\vspace{.10in}
\section{Introduction}\label{section int}
The system of differential equations
 \begin{equation}
 \label{darboux}
\left \{ \begin{array}{l}
\frac{dt_1}{dz}+\frac{dt_2}{dz}=t_1t_2\\
\frac{dt_2}{dz}+\frac{dt_3}{dz}=t_2t_3\\
\frac{dt_1}{dz}+\frac{dt_3}{dz}=t_1t_3
\end{array} \right.,
\end{equation}
appeared in 1878 in the work of Gaston Darboux \cite{da78},
 where he was treating on the  curvilinear coordinates
and orthogonal systems. The problem that he was
working on it is as follow: \emph{Let $A$ and $B$ be two
fixed surfaces in 3-dimensional Euclidean space
$\mathbb{R}^3$. Suppose that $\Sigma$ is a family of
surfaces parallel to $A$, and
$\Sigma '$ is another family  of surfaces parallel to $B$. Is there a third family of surfaces
parameterized by $u$ such that intersects $\Sigma$ and $\Sigma
'$ orthogonally?}
The more interesting case of this problem is when the family
$(u)$ is of the second degree and Darboux proved that in this case
this family is given by
\[
\frac{x_1^2}{t_1(u)}+\frac{x_2^2}{t_2(u)}+\frac{x_3^2}{t_3(u)}=1,
\]
in which $x_1,x_2,x_3$ are coordinates of $\R^3$, and $t_1,t_2,t_3$ are functions of $u$ given by the
following equation
\begin{equation}\label{eq darboux}
t_3(\frac{dt_1}{du}+\frac{dt_2}{du})=t_2(\frac{dt_1}{du}+\frac{dt_3}{du})=t_1(\frac{dt_2}{du}+\frac{dt_3}{du}).
\end{equation}
Therefore, the system of equations (\ref{darboux}) is a particular
case of the equation (\ref{eq darboux}).

In 1881, G. Halphen \cite{ha81} studied the system of
differential equations (\ref{darboux}) in $\C^3$. He
proved that this system satisfies an important invariant
property. To express this invariant property, for the constants
$a,b,a',b'$, let
 \begin{equation} \label{eq change var1}
w=\frac{az+b}{a'z+b'}\qquad \& \qquad t_i=-\frac{2a'}{a'z+b'}+\frac{ab'-ba'}{(a'z+b')^2}s_i,\,\,\, i=1,2,3.
\end{equation}
By substituting (\ref{eq change var1}) in the system (\ref{darboux}), we have
\begin{equation}
 \label{darboux1}
\left \{ \begin{array}{l}
\frac{ds_1}{dw}+\frac{ds_2}{dw}=s_1s_2\\
\frac{ds_2}{dw}+\frac{ds_3}{dw}=s_2s_3\\
\frac{ds_1}{dw}+\frac{ds_3}{dw}=s_1s_3
\end{array} \right.,
\end{equation}
from which it follows that the system (\ref{darboux}) is invariant under the change
of variables (\ref{eq change var1}).
Therefore, to find a general solution of (\ref{darboux}), it is enough
to apply (\ref{eq change var1}) to  a particular solution of (\ref{darboux1}). Halphen gave a
solution of the system (\ref{darboux}) in terms of the
logarithmic derivatives of the null theta functions; namely
\begin{eqnarray*}\small
t_1&=& 2(\ln \theta_4(0|z))',\\
t_2&=&2(\ln \theta_2(0|z))', \ \ \ \ \  '=\frac{\partial}{\partial z},\\
t_3&=&2(\ln \theta_3(0|z))'.
\end{eqnarray*}
where
$$\small
\left \{ \begin{array}{l}
\theta_2(0|z):=\sum_{n=-\infty}^\infty q^{\frac{1}{2}(n+\frac{1}{2})^2}
\\
\theta_3(0|z):=\sum_{n=-\infty}^\infty q^{\frac{1}{2}n^2}
\\
\theta_4(0|z):=\sum_{n=-\infty}^\infty (-1)^nq^{\frac{1}{2}n^2}
\end{array} \right.,
\ q=e^{2\pi i z},\  \textrm{Im}(z)>0.
$$

F. Brioschi \cite{brio} in 1881 studied the following
extension of the system (\ref{darboux})
\begin{equation}
 \label{brio}
\left \{ \begin{array}{l}
\frac{dt_1}{dz}+\frac{dt_2}{dz}=t_1t_2+\varphi(z)\\
\frac{dt_2}{dz}+\frac{dt_3}{dz}=t_2t_3+\varphi(z)\\
\frac{dt_1}{dz}+\frac{dt_3}{dz}=t_1t_3+\varphi(z)
\end{array} \right.,
\end{equation}
in which $\varphi(z)$ is a function of $z$. Again in 1881,
Halphen in \cite{ha81-1} introduced and investigated a class of
differential equations with (\ref{brio}) as a particular case. In the case of three variables, he showed that this
class is given by
\begin{equation}
 \label{eq hal}
\left \{ \begin{array}{l}
\frac{dt_1}{dz}=a_1t_1^2+(\lambda-a_1)(t_1t_2+t_1t_3-t_2t_3)\\
\frac{dt_2}{dz}=a_2t_2^2+(\lambda-a_2)(t_2t_3+t_2t_1-t_3t_1)\\
\frac{dt_3}{dz}=a_3t_3^2+(\lambda-a_3)(t_3t_1+t_3t_2-t_1t_2)
\end{array} \right.,
\end{equation}
where, $a_1,a_2,a_3,\lambda$ are constants. He
proved that the system (\ref{eq hal}) also satisfies the
invariant property and it is in a direct relationship with the
Gauss hypergeometric equation (see \cite{ha81-0}). One can see that
the system (\ref{eq hal}) is equivalent to the system
(\ref{darboux}), when $a_1=a_2=a_3=0$ and $\lambda=1$. If we look at the system (\ref{eq hal}) as a vector field in $\C^3$, then it is a semi-complete vector field. In this context, an extension of Halphen vector field, namely Halphen type vector field, was introduced by Adolfo Guillot in \cite{ado1,ado2}.

In 1916 Ramanujan \cite{ra16} introduced another system of differential equations as follow
\begin{equation} \label{ramanujan}
\textrm{R}:\, \left \{ \begin{array}{l}
\frac{dr_1}{d\tau}=r_1^2-\frac{1}{12}r_2 \\
\frac{dr_2}{d\tau}=4r_1r_2-6r_3 \\
\frac{dr_3}{d\tau}=6r_1r_3-\frac{1}{3}r_2^2
\end{array} \right.,
\end{equation}
that is in a close relationship with Darboux-Halphen differential equation (\ref{darboux}).
He verified that the Eisenstein series $\frac{2\pi i}{12}E_2(\tau),$ $12(\frac{2\pi i}{12})^2E_4(\tau),8(\frac{2\pi i}{12})^3E_6(\tau)$ satisfy (\ref{ramanujan}), where
\begin{align}
&E_{2j}(\tau):=1-\frac{4j}{B_{2j}}\sum_{r=1}^\infty \sigma_{2j-1}(r)q^r,\,\,\,\, q=e^{2\pi i\tau}, \nonumber \\
                                             & \sigma_{i}(n):=\sum_{d|n}d^i  \nonumber,
\end{align}
and $B_k$'s are Bernoulli's numbers.
In (\ref{eq morphism}) we will see a relationship between the systems of equations
(\ref{darboux}) and (\ref{ramanujan}).\\

\emph{Calabi-Yau manifolds} are defined as compact connected
K\"{a}hler manifolds whose canonical bundle is trivial, though
many other equivalent definitions are sometimes used. They were named "\cy manifold"
by Candelas et al. (1985) \cite{canetal} after E. Calabi (1954)
\cite{calabi1, calabi2},
who first studied them, and S. T. Yau (1976) \cite{yau1}, who proved
the Calabi conjecture that says \cy manifolds accept Ricci flat metrics.
In this text, we suppose that for an $n$-dimensional \cy manifold $\h^{p,0}=0$, $0<p<n$, where
$\h^{p,q}$ refers to $(p,q)$-th Hodge number of \cy manifold (see \S \ref{section pfe}).
It is clear that the connectedness of a \cy manifold and the triviality of its canonical bundle imply that $\h^{0,0}=\h^{n,0}=1$.
In order to explain the generalization of Darboux-Halphen-Ramanujan vector fields, $\DHR$ for short, we consider the family of one-dimensional Calabi-Yau manifolds, which are elliptic curves, and for more details the reader refers to \cite{ho06-1}.

Let $E$ be an elliptic curve over $\C$. Then the Hodge
filtration $F^\bullet H^1$ of
the first de Rham cohomology group $H^1_\dR(E)$ is given as follow,
\[
\{0\}=F^{2}\subset F^1\subset F^0=H^{1}_\dR(E),\ \ \dim F^i=2-i,
\]
where $F^1\subset H^{1}_\dR(E)$ includes classes of holomorphic closed $1$-forms on $E$. Let \Tm be the moduli of the pair $(E,[\alpha_1,\alpha_2])$, in
which $\alpha_1\in F^1$,  $\alpha_2\in F^0\setminus F^1$, and
the intersection form matrix in $\alpha_i$'s is as follow
\[\small
\left( \langle \alpha_i,\alpha_j \rangle \right)_{1\leq i,j\leq 2}=\left(\begin{array}{cc}
                                                                     0 & 1 \\
                                                                     -1 & 0
                                                                   \end{array}\right).
\]
To be more precise, let $E_0:=E\setminus \{\infty\}$ be the affine curve that its Weierstrass presentation is given as follow
\[
E_0=\{(x,y)\in \C^2 |\, f(x,y):=y^2-4(x-t_1)^3+t_2(x-t_1)+t_3=0\}.
\]
Then $\alpha_1$ and $\alpha_2$, resp., are
induced by $[\frac{dx}{y}]$ and $[\frac{xdx}{y}]$, resp., where
$[\frac{dx}{y}]$ and $[\frac{xdx}{y}]$ are generators of the
first de Rham cohomology $H^1_\dR(E_0)$ of the affine curve
$E_0$. Since $H^1_\dR(E)\cong
H^1_\dR(E_0)$,  it follows that $\alpha_1$ and $\alpha_2$ are generators of
$H^1_\dR(E)$; and hence $\alpha_1\wedge \alpha_2\neq 0$. It is seen that \Tm is a
3-dimensional space, and there exist a unique vector field $\H$
on \Tm such that the composition of Gauss-Manin connection (see \S \ref{subsection gm})
$$\nabla:H_\dR^1(E/\Tm)\to \Omega^1_\Tm\otimes_{\O_\Tm}H_\dR^1(E/\Tm),$$
with $\H$ satisfies the following:
\[\small
\nabla_{\H}\left( \begin{array}{c}
                     \alpha_1 \\
                     \alpha_2
                   \end{array}
 \right)=\left(\begin{array}{cc}
                                                                     0 & -1 \\
                                                                     0 & 0
                                                                   \end{array}\right)\left( \begin{array}{c}
                     \alpha_1 \\
                     \alpha_2
                   \end{array}
 \right).
\]
Roughly speaking, for $y^2=4(x-t_1)(x-t_2)(x-t_3)$ we have $\Tm=\Tm_\textrm{DH}:=\{(t_1,t_2,t_3)\in \mathbb{C}^3| t_1\neq t_2\neq t_3\}$ and $\H=\textrm{DH}$ is
given by the following system
\begin{equation}\label{darboux2}
\textrm{DH}:\,\left \{ \begin{array}{l}
\frac{dt_1}{dz}=t_1(t_2+t_3)-t_2t_3\\
\frac{dt_2}{dz}= t_2(t_1+t_3)-t_1t_3 \\
\frac{dt_3}{dz}= t_3(t_1+t_2)-t_1t_2
\end{array} \right.,
\end{equation}
which is an special case of the system (\ref{eq hal}) introduced by Darboux-Halphen with $a_1=a_2=a_3=0$ and $\lambda=1$. Or
equivalently for $y^2=4(x-r_1)^3+t_2(x-r_1)+r_3$ we obtain $\Tm=\Tm_\textrm{R}:=\{(r_1,r_2,r_3)\in\mathbb{C}^3|27r_3^2-r_2^3=0\},$ and $\H=\textrm{R}$ is presented by system (\ref{ramanujan}) introduced by Ramanujan (for details see \cite[Proposition~3.8]{ho06-1}). The
algebraic morphism $\phi:\Tm_{\rm DH}\to \Tm_{\rm R}$ defined by
\begin{equation}\label{eq morphism}
\phi: (t_1,t_2,t_3)\mapsto (T, 4\sum_{1\leq i<j\leq 3}(T-t_i)(T-t_j), 4(T-t_1)(T-t_2)(T-t_3)),
\end{equation}
where $\ T:=\frac{1}{3}(t_1+t_2+t_3),$
connects two systems (\ref{ramanujan}) and (\ref{darboux2}), i.e., $\phi_\ast {\rm DH}={\rm
R}$.\\

After these works, H. Movasati \cite{ho21} considered a
one parameter family of \cy 3-folds, which is known as the
family of mirror quintic 3-folds, and studied on it. If $W$ is
a mirror quintic 3-fold, then the Hodge filtration of
$H^3_\dR(W)$ is given as follow
\[
\{0\}=F^{4}\subset F^3\subset\ldots \subset F^0=H^{3}_\dR(W),\ \ \dim F^i=4-i.
\]
The complex moduli of $W$ is one dimensional that we parameterize it by $z$. There is a nowhere vanishing holomorphic 3-form $\omega\in F^3$ such that
the Picard-Fuchs equation associated with it is given by
\begin{align}
\L&=\vartheta^4-5^5z(\vartheta+\frac{1}{5})(\vartheta+\frac{2}{5})(\vartheta+\frac{3}{5})(\vartheta+\frac{4}{5})\nonumber,
\end{align}
in which $\vartheta:=\nabla_{z\frac{\partial}{\partial z}}$ is the composition of Gauss-Manin connection $\nabla$ with the vector field $z\frac{\partial}{\partial z}$.
Movasati treated on the moduli space \Tm of
the pair $(W,[\alpha_1,\alpha_2,\alpha_3,\alpha_4])$, where
$\alpha_i\in F^{4-i}\setminus F^{5-i}$, and the intersection form matrix in $\alpha_i$'s is given as follow
\[\small
\left( {\langle \alpha_i,\alpha_j\rangle}\right)_{1\leq i,j\leq 4}=\begin{pmatrix}
 0&0&0&1\\
0&0&1&0\\
0&-1&0&0\\
-1&0&0&0
\end{pmatrix}.
\]
He proved that \Tm is a 7-dimensional space and there is a
unique vector field $\H$ and a unique meromorphic function $y$
on \Tm such that,
\[\small
\nabla _{\H} \left( \begin{array}{c}
  \alpha_1 \\
  \alpha_2 \\
  \alpha_3 \\
  \alpha_4
\end{array} \right) = \left( {\begin{array}{*{20}c}
   0 & 1 & 0 & 0  \\
   0 & 0 & y & 0  \\
   0 & 0 & 0 & { - 1}  \\
   0 & 0 & 0 & 0  \\
\end{array}} \right) \left( \begin{array}{c}
  \alpha_1 \\
  \alpha_2 \\
  \alpha_3 \\
  \alpha_4
\end{array} \right).
\]
Indeed he expressed $\H$ and $y$ explicitly, and he showed that $y$ is related with the normalized
Yukawa coupling, whose was introduced by Candelas et
al. in \cite{can91}. They computed  the coefficients of the
$q$-expansion of the  normalized  Yukawa coupling for quintic 3-folds in
$\P^4$, that are conjectured to be the Gromov-Witten invariants of rational
curves on a quintic $3$-fold in $\P^4$. \\

After what we saw about the family of \cy 1-folds and the
family of mirror quintic 3-folds, it is natural to ask whether there
exist such a moduli space \Tm and such a vector field $\H$ in higher dimensions. In the
present paper we give a positive answer to this question. To do this, we fix an $n$-dimensional \cy manifold
$W$. We suppose that the complex deformation of $W$ is
given by a one parameter family $\pi:\W\to P$, where $P$ is a one-dimensional quasi-projective variety
parameterized by $z$. Moreover, we assume that the $n$-th relative de Rham cohomology group $H^n_\dR(\W/P)$ is
$n+1$-dimensional and the Picard-Fuchs equation $\L$ associated with the unique nowhere vanishing
holomorphic $n$-form $\omega\in \F^n$ is given by
\begin{equation}
\L=\vartheta^{n+1}-a_n(z)\vartheta^n-\ldots-a_1(z)\vartheta-a_0(z),
\end{equation}
where $\F^\bullet H^n$ is
the Hodge filtration of $H^n_\dR(\W/P)$, $\vartheta:=\nabla_{z\frac{\partial}{\partial z}}$ and
$a_i(z)\in \Q(z)$, $i=0,1,\ldots,n$. We provide the first result in the following proposition.

\begin{prop}\label{prop 1}
The Picard-Fuchs equation $\L$ is self-dual.
\end{prop}

Before stating the main
theorem of this paper, we fix the $(n+1)\times (n+1)$ matrix
$\Phi$ as follow. If $n$ is an \emph{odd} integer, then set
\begin{equation}\label{eq phi}\small
\Phi:=\left( {\begin{array}{*{20}c}
   {0_{\frac{n+1}{2}} } & {J_{\frac{n+1}{2}}}   \\
   { - J_{\frac{n+1}{2}}}  & {0_{\frac{n+1}{2}}}   \\
\end{array}} \right),
\end{equation}
where for $k\in \mathbb{N}$, $0_{k}$ denotes a $k\times
k$ block of zeros, and $J_k$ is the following $k\times k$ block
\begin{equation}\small
J_{k }  := \left( {\begin{array}{*{20}c}
   0 & 0 &  \ldots  & 0 & 1  \\
   0 & 0 &  \ldots  & 1 & 0  \\
    \vdots  &  \vdots  &  {\mathinner{\mkern2mu\raise1pt\hbox{.}\mkern2mu
 \raise4pt\hbox{.}\mkern2mu\raise7pt\hbox{.}\mkern1mu}}  &  \vdots  &  \vdots   \\
   0 & 1 &  \ldots  & 0 & 0  \\
   1 & 0 &  \ldots  & 0 & 0  \\
\end{array}} \right).
\end{equation}
If $n$ is an \emph{even} integer, then $\Phi:=J_{n+1}$. Also we suppose that $F^\bullet H^n$ denotes the Hodge filtration of $H^n_\dR(W;\C)$ given as follow:
\[
F^\bullet H^n:\,\ \{0\}=F^{n+1}\subset F^n\subset \ldots \subset F^1\subset F^0=H^{n}_\dR(W;\C).
\]

\begin{theo}\label{theo 1 int}
Let $W$ be the \cy $n$-fold given above and \Tm be the moduli
of $(W,[\alpha_1,\alpha_2,\ldots,\alpha_{n+1}])$, where $\{\alpha_i\}_{i=1}^{n+1}$ is a basis of $H^{n}_\dR(W;\C)$ satisfying
\begin{equation} \label{eq 3714 2}
\alpha_i\in F^{n+1-i}\setminus F^{n+2-i},\,\, i=1,2,\ldots,n+1,
\end{equation}
and the intersection form matrix in $\alpha_i$'s is subject to the condition:
\begin{equation}\label{eq phi n}
\left( \langle\alpha_i,\alpha_j\rangle \right)_{1\leq i,j\leq n+1}=\Phi.
\end{equation}
Then there exist a unique vector field $\H$ and unique
meromorphic functions $y_i$, $i=1,2,\ldots,n-2$, on \Tm such
that the composition of Gauss-Manin connection $\nabla$ with
the vector field $\, \H\, $ satisfies:
\begin{equation}\label{eq gm halphenn}
\nabla_\H\alpha=Y\alpha,
\end{equation}
in which
\[
\alpha=\left( {\begin{array}{*{20}c}
   {\alpha _1 } & {\alpha _2 } & {\ldots } & {\alpha _{n+1} }  \\
\end{array}} \right)^\tr,
\]
and
\begin{equation}\label{eq y n}\small
Y=\left(
         \begin{array}{cccccc}
           0 & 1 & 0 & \ldots & 0 & 0 \\
           0 & 0 & y_1 & \ldots & 0 & 0 \\
           \vdots & \vdots & \vdots & \ddots & \vdots & \vdots \\
           0 & 0 & 0 & \ldots & y_{n-2} & 0 \\
           0 & 0 & 0 & \ldots & 0 & -1 \\
           0 & 0 & 0 & \ldots & 0 & 0 \\
         \end{array}
       \right).
\end{equation}
Moreover we have,
{\rm $$ \dim \Tm=\left \{
\begin{array}{l}
\frac{(n+1)(n+3)}{4}+1;\,\, \quad  \textrm{\rm if \textit{n} is odd} \\
\frac{n(n+2)}{4}+1;\,\,\,\,\quad\quad \textrm{\rm if \textit{n} is even}
\end{array} \right..
$$}
\end{theo}

As we saw above, the system of ordinary differential equations given by $\H$ is an extension of systems of differential equations introduced by Darboux, Halphen and Ramanujan.

\begin{defi}
{\rm The vector field $\H$ introduced in Theorem \ref{theo 1 int}, is called \emph{Darboux-Halphen-Ramanujan}, $\DHR$ for brevity, \emph{vector field}.}
\end{defi}

The structure of this article is prepared as follow. First, in Section \ref{section do} we give an algorithm to find the existence relationships among
coefficients of a self-dual linear differential equation of an arbitrary
degree. In particular we provide these relationships in degrees three and five. Section \ref{section pfe} contains a brief summary of some basic facts. After fixing some notations and assumptions, the proof of Proposition \ref{prop 1} is given in \S \ref{subsection sd}. Finally, Section \ref{section dhr} is devoted to the proof of Theorem \ref{theo 1 int}. In this section the proof is divided to the even case and odd case depending to the dimension of \cy manifold. And also we present \DHR vector field explicitly in dimensions three and five.

\begin{rem}
As we will see in \S\ref{section dhr}, to prove Theorem \ref{theo 1 int} we introduce several matrices and matrix equations. Recently we discovered that they are in a close relationship with Birkhoff factorization given in quantum cohomology (see {\rm \cite{guest}}). In fact, if we talk in physicists language, our work is in B-model and the Birkhoff factorization is discussed in A-model and mirror symmetry gives the existence relationships between them.
\end{rem}

{\bf Acknowledgment.} Here I would like to express my very great appreciation to Hossein Movasati, my Ph.D. supervisor, who always was available and I used his valuable and constructive suggestions and helps during the planning and development of this  work. I wish to thank IMPA for preparing such an excellent academic environment. This work has been done during my Ph.D. and I am grateful to have economic supports of "CNPq-TWAS Fellowships Programme" during this period.

\section{Self-Dual Linear Differential Equation}\label{section do}
In this section by $R$ we mean the simple commutative
differential ring $\C[z]$, with quotient field $k:=\C(z)$ and derivative
$(.)'$; and $R[\partial]$ is the ring of differential
operators where $\partial$ is the usual derivation
$\frac{\partial}{\partial z}$ or logarithmic derivation
$z\frac{\partial}{\partial z}$. It is not difficult to check that
$k[\frac{\partial}{\partial z}]$ and
$k[z\frac{\partial}{\partial z}]$ are isomorphic, hence we can
freely switch between these two differential rings. The pair $(M,\partial)$ refers
to a differential $R$-module, i.e., $M$ is a finitely generated $R$-module and $\partial: M
\to M$ is a map satisfying $\partial(m + n) = \partial(m) +\partial(n)$ for every $m, n \in M$; and
$\partial(fm) = f'm + f\partial(m)$ for every $f\in R$ and every $m\in M$. For more details the reader can see \cite{pusi}.

\begin{defi} \rm
Let $(M,\partial)$ be a differential $k$-module. Then for each
$m\in M$ we define the \emph{evaluation map} $ev_m:
k[\partial]\to M$ by $\sum_{i = 0}^n a_i\partial^i\mapsto \sum_{i = 0}^na_i\partial^im$.
The monic generator of the kernel of $ev_m$ as a left ideal is
called the \emph{minimal operator} of $m$ over $k[\partial]$.
Furthermore, we call $m$ a \emph{cyclic vector} of $M$ if the
degree of its minimal operator equals the $k$-dimension of $M$,
i.e. the set $\{m, \partial m, . . . ,
\partial^{\dim_k(M)-1}m\}$ is a $k$-basis of M. We call a pair
$(M, e)$ consisting of a differential module $M$ and a cyclic
vector $e\in M$ a \emph{marked differential module}.
\end{defi}

By a result due to N. Katz (see \cite[\S~2.1]{pusi}), there is a one to one
correspondence between monic differential equations $L\in
k[\partial]$ and marked differential modules $(M,e)$. More precisely,
each differential $k$-module $M$ has a cyclic vector, and in
particular there is a differential equation $L\in k[\partial]$
such that $M$ is isomorphic to $k[\partial]/k[\partial]L$. Thus we
can assume
$
L =\partial^{n+1} +\sum_{i = 0}^n a_i\partial^i\in \mathbb{Q}(z)[\partial],
$
is an irreducible monic differential equation and
\begin{equation}
(M_L, e)\cong(\mathbb{C}(z)[\partial]/\mathbb{C}(z)[\partial]L , [1]),
\end{equation}
is its corresponding marked differential
$\mathbb{C}(z)$-module. The \emph{dual equation} $\check{L}$ of $L$ is defined as follow
\begin{equation}
\check{L}=\sum\limits_{i = 0}^{n+1} {(-1)^{n-i}\partial ^i a_i },\,\,\ a_{n+1}=1.
\end{equation}

\begin{defi}\label{defi property p}\rm
It is said that $L$ satisfies \emph{property} (P), if there is
a non-degenerate form $\langle .,. \rangle:M_L\times M_L\to
\mathbb{C}(z)$ such that
\begin{description}
  \item[(i)] $\langle .,.\rangle$ is a $(-1)^n$-symmetric
      form, i.e. $\langle .,. \rangle \in  {\rm
      Hom_{\mathbb{C}(z)[\partial]}(Sym^2}M_L,\mathbb{C}(z))$
      if $n$ is even, and $\langle .,. \rangle \in  {\rm
      Hom_{\mathbb{C}(z)[\partial]}(\bigwedge^2}M_L,\mathbb{C}(z))$
      if $n$ is odd.
  \item[(ii)] $\langle e, \partial^ie\rangle=0$ for $i=0,1,
      \ldots, n-1.$
\end{description}
\end{defi}

We state a proposition that gives an equivalence
condition for property (P) and for a proof see \cite{bogn}.
Note that for $\psi\in \C(z)$, the operator $\partial \psi$ is
given as $\partial \psi=\partial(\psi)+\psi\partial$,
and for convenient we denote by $\psi'=\partial(\psi)$, so
$\psi^{(i)}=\underbrace{\partial(\partial(\ldots(\partial }_{i
- {\rm times}}(\psi))\ldots))$.

\begin{prop} \label{prop bog}
The equation $L$ satisfies the property {\rm (P) }if and only
if $L$ is self-dual, i.e., there is an $0\neq \psi\in
\mathbb{C}(z)$, such that
\begin{equation}
L\psi=\psi \check{L}.
\end{equation}
\end{prop}

Using Proposition \ref{prop bog},  we give an
algorithm to find the relationships that exist among
coefficients $a_i$'s. Let $L=\sum_{i=0}^{n+1}a_i\partial^i$, with
$a_{n+1}=1$, be a linear differential equation satisfying
property {\rm (P)}. Suppose that $n=2m \, {\rm or}\, 2m+1$, for
a positive integer $m$. Then coefficients
$a_{n-2},a_{n-4},\ldots,a_{n-2m}$ depend to the rest of the
coefficients and their derivations. First using the induction, one can
easily verify that for $\psi \in \C(z)$
\[\small
\partial^j\psi=\sum\limits_{i=0}^j\left(
                                                               \begin{array}{c}
                                                                 j \\
                                                                 i \\
                                                               \end{array}
                                                             \right)
                                                             \psi^{(j-i)}\partial^i.
\]
Therefore, it follows that
\begin{equation}\label{eq lcheck}\small
  \check{L} =\sum\limits_{i=0}^{n+1}\left(\sum\limits_{j=i}^{n+1}(-1)^{n+1-j}\left(
                                                               \begin{array}{c}
                                                                 j \\
                                                                 i \\
                                                               \end{array}
                                                             \right)
  a_j^{(j-i)}\right)\partial^i ,
\end{equation}
and
\begin{equation}\label{eq lpure}\small
  L\psi =\sum\limits_{i=0}^{n+1}\left(\sum\limits_{j=i}^{n+1}\left(
                                                               \begin{array}{c}
                                                                 j \\
                                                                 i \\
                                                               \end{array}
                                                             \right)
  a_j\psi^{(j-i)}\right)\partial^i .
\end{equation}
If we substitute (\ref{eq lcheck}) and (\ref{eq lpure}) in $L\psi=\psi\check{L}$, then we
have
\begin{equation}\label{eq 42214 1}\small
\sum\limits_{i=0}^{n+1}\left(\sum\limits_{j=i}^{n+1}\left(
                                                               \begin{array}{c}
                                                                 j \\
                                                                 i \\
                                                               \end{array}
                                                             \right)
  a_j\psi^{(j-i)}\right)\partial^i=
\sum\limits_{i=0}^{n+1}\psi \left(\sum\limits_{j=i}^{n+1}(-1)^{n+1-j}\left(
                                                               \begin{array}{c}
                                                                 j \\
                                                                 i \\
                                                               \end{array}
                                                             \right)
  a_j^{(j-i)}\right)\partial^i.
\end{equation}
Now by comparing the coefficient of $\partial^{n}$ in (\ref{eq
42214 1}), we express $\psi'$ and $\psi^{(i)}$'s in terms of
$\psi$, $a_n$ and derivations of $a_n$ as follows
\begin{align}\small
&\psi'=-\frac{2}{n+1}a_n\psi, \nonumber \\
&\psi''=\left((-\frac{2}{n+1})^2a_n^2-\frac{2}{n+1}a'_n\right)\psi,  \\
&\psi'''=\left( (-\frac{2}{n+1})^3a_n^3+3(-\frac{2}{n+1})^2a_na'_n -\frac{2}{n+1}a''_n\right)\psi, \nonumber \\
&\vdots \nonumber
\end{align}
and we substitute $\psi^{(i)}$' in the left side of (\ref{eq 42214 1}).
In order to state $a_{n-2k},\, k=1,2,\ldots,m$, as an equation
of $a_n,a_{n-1},a_{n-3},\ldots,a_{n-(2k-1)}$ and their
derivations it is enough to compare the coefficient of
$\partial^{n-2k}$ of both sides of (\ref{eq 42214 1}), which yields
\[\small
\sum\limits_{j={n-2k}}^{n+1}\left(
                                                               \begin{array}{c}
                                                                 j \\
                                                                 {n-2k} \\
                                                               \end{array}
                                                             \right)
  a_j\psi^{(j-({n-2k}))}=\left(\sum\limits_{j={n-2k}}^{n+1}(-1)^{n+1-j}\left(
                                                               \begin{array}{c}
                                                                 j \\
                                                                 {n-2k} \\
                                                               \end{array}
                                                             \right)
  a_j^{(j-({n-2k}))}\right)\psi.
\]
Therefore,
\begin{align}\small
2a_{n-2k}=\sum\limits_{j={n-2k}+1}^{n}(-1)^{n+1-j} & \left(
                                                               \begin{array}{c}
                                                                 j \\
                                                                 {n-2k} \\
                                                               \end{array}
                                                             \right)
  a_j^{(j-({n-2k}))}  \nonumber \\
  - \sum\limits_{j={n-2k+1}}^{n+1} & \left(
                                                               \begin{array}{c}
                                                                 j \\
                                                                 {n-2k} \\
                                                               \end{array}
                                                             \right)
  a_j(\psi^{(j-({n-2k}))}/\psi). \nonumber
\end{align}
For example if $k=1$, then $a_{n-2}$ is given as follow
\[\small
a_{n-2}=\frac{n-1}{n+1}a_{n-1}a_n-\frac{n(n-1)}{2(n+1)}a_na'_n-\frac{n(n-1)}{3(n+1)^2}a_n^3+\frac{(n-1)}{2}a'_{n-1}-\frac{1}{12}n(n-1)a''_n.
\]
As a result of this algorithm we provide the following lemma.

\begin{lemm} \label{lemm pp}
Let $L=\sum_{i=0}^{n+1}a_i\partial^i$, with $a_{n+1}=1$, be a
linear differential equation satisfying property {\rm (P)}.
Then followings hold:
\begin{description}
  \item[(i)] If $n=3$, then
  \[\small
  a_1=\frac{1}{2}a_2a_3-\frac{3}{4}a_3a_3'-\frac{1}{8}a_3^3+a_2'-\frac{1}{2} a_3''.
  \]
  \item[(ii)] If $n=5$, then
  \begin{align}\small
a_3&=\frac{2}{3}a_4a_5-\frac{5}{3}a_5a_5'-\frac{5}{27}a_5^3+2a_4'-\frac{5}{3}a_5'', \nonumber \\
a_1&=a_2'-a_4'''+a_5^{(4)}-a_4^{(2)}a_5- a_4' a_5'+\frac{5}{3}a_5(a_5')^2+\frac{1}{3}a_2a_5 \nonumber \\
&-\frac{1}{27}a_4a_5^3
+\frac{10}{27}a_5^3a_5'+\frac{1}{81}a_5^5-\frac{1}{3} a_4'a_5^2-\frac{1}{3}a_4a_5a_5'+\frac{10}{9}a_5^2a_5'' \nonumber   \\
&+\frac{10}{3} a_5'a_5''-\frac{1}{3}a_4a_5''+\frac{5}{3}a_5a_5''' \nonumber.
\end{align}
\end{description}
\end{lemm}

\section{Picard-Fuchs Equation as a Self-Dual Linear Differential Equation}\label{section pfe}

In this section $\pi:\W\to P$ refers to a family of $n$-dimensional compact K\"ahler manifolds, i.e., $\pi$ is a holomorphic proper submersion of complex manifolds $\W$ and $P$ such that for any $z\in P$, $W_z:=\pi^{-1}(z)$ is an $n$-dimensional compact K\"ahler manifold. If we denote the $k$-th de Rham cohomology group of $W_z$ by $H^k_\dR(W_z)$, then \emph{de Rham Lemma} gives the isomorphism $H^k_{dR}(W_z)\cong H^k(W_z,\R)$, or equivalently $H^k_{dR}(W_z;\C)\cong H^k(W_z,\C)$ where $H^k_{dR}(W_z;\C)$ denotes the complexified de Rham cohomology group. Here $\b_k(W_z):=\dim H^k_{dR}(W_z;\C)$ stands for the $k$-th betti number of $W_z$. Also by \emph{Hodge decomposition theorem} we have
\begin{equation}\label{eq hdt}
H^k_{dR}(W_z;\C)=\bigoplus\limits_{p+q=k}H^{p,q}(W_z),
\end{equation}
in which $H^{p,q}(W_z)$ is $(p,q)$-th Dolbeault cohomology and by \emph{Dolbeault's
theorem} we have the isomorphism $H^{p,q}_{}(W_z)\cong H^q(W_z,\Omega^p_{W_z})$. We denote by $\h^{p,q}(W_z):=\dim H^{p,q}(W_z)$, that is called \emph{$(p,q)$-th Hodge number} of $W_z$. By defining $$F^p(W_z):=\bigoplus\limits_{p\leq r\leq n}H^{r,n-r}(W_z),\,\, 0\leq p\leq n,$$ we yield the following decreasing filtration which is known as the \emph{Hodge filtration} of $H^k_\dR(W_z)$,
\begin{equation}\label{eq d hodge f1}
F^\bullet H^k(W_z):\,\,\{0\}=F^{k+1}(W_z)\subset F^k(W_z)\subset\ldots\subset F^0(W_z)=H^k_\dR(W_z;\C).
\end{equation}
We can consider the family $\W$ as a complex deformation of $W:=W_0$, $0\in P$. As one can find in standard texts of complex geometry, e.g. \cite{vo02}, up to replacing $P$ by a neighborhood of the base point $0$,  $\b_k(W_z)=\b_k(W)$ and $\h^{p,q}(W_z)=\h^{p,q}(W)$ for any $z\in P$. Hence simply we can write $\b_k$ and $\h^{p,q}$ instead of $\b_k(W_z)$ and $\h^{p,q}(W_z)$. Also one can see that $\b_k=\sum_{p+q=k}\h^{p,q}$, $\h^{p,q}=\h^{q,p}$ and $\h^{p,q}=\h^{n-q,n-p}$.

\subsection{Gauss-Manin Connection and Griffiths Transversality}\label{subsection gm}

Consider the sheaf $R^k\pi_\ast\underline{\C_{\W}}$ on $P$, where $\underline{\C_{\cal W}}$ is the constant sheaf on $\W$ with fibers $\C$ and
$R^k\pi_\ast$ refers to $k$-th derived functor of the pushforward. For any $z\in P$, we have the following presentation of the stalks of $R^k\pi_\ast\underline{\C_{\W}}$
\[
(R^k\pi_\ast\underline{\C_{\W}})_z \simeq H^k(W_z,\C)\mathop  \simeq  H_\dR^k(W_z;\C).
\]
Hence $R^k\pi_\ast\underline{\C_{\W}}$ is a locally
constant sheaf on $P$. Formally speaking,
$R^k\pi_\ast\underline{\C_{\W}}$ is the sheaf associated to
the presheaf $U\mapsto H^k(\pi^{-1}(U),\C)$ (see \cite{grjo}). In fact, for a
contractible open subset $U\subset P$, by Ehresmann Lemma
$\pi^{-1}(U)\cong U\times W_z$ for some $z\in P$, so
$H^k(\pi^{-1}(U),\C)\simeq H^k(W_z,\C)$. By defining
\begin{equation} \label{eq rdrcg}
H^k_\dR({\W}/P):=R^k\pi_\ast\underline{\C_{\W}}\otimes_\C \O_P,
\end{equation}
which is a holomorphic vector bundle on $P$, then for any $z\in P$, $H_\dR^k({\W}/P)_z\cong H_\dR^k(W_z;\C)$.
\begin{defi}
{\rm The holomorphic vector bundle $H^k_\dR({\W}/P)$
defined in (\ref{eq rdrcg}) is called \emph{$k$-th relative de
Rham cohomology group}. The unique integrable connection
\[
\nabla^\GM:H^k_{\dR}({\W}/P)\to \Omega_P^1 \otimes_{\O_P} H^k_{\dR}({\W}/P),
\]
whose flat sections coincides with  $R^k\pi_\ast\underline{\C_{\W}}$ is known as \emph{Gauss-Manin connection}. }
\end{defi}
For a vector field $v$ on $P$, consider the map $v\otimes {\rm Id}: \Omega^1_P\otimes_{\O_P} H_\dR^k({\W}/P)\to
 H_\dR^k(\W/P)$. Then by composing the Gauss-Manin connection $\nabla^\GM$ with
$v\otimes {\rm Id}$ we define
\begin{align}\label{eq gmovf}
\nabla^\GM_v:& H_\dR^k(\W/P)\to H_\dR^k(\W/P) \\
&\nabla^\GM_v:=(v\otimes {\rm Id})\circ \nabla^\GM.\nonumber
\end{align}
From now on, if no confusion arises, we denote the
Gauss-Manin connection by $\nabla$ instead of  $\nabla^\GM$.
\begin{rem}\label{obs gm mat}
{\rm The $k$-th relative de Rham cohomology group $H^k_\dR(\W/P)$
is locally free of finite rank, say $m$. Let $\{\omega_j\}_{j=1}^m$ be a local
 frame  of $H^k_\dR(\W/P)$ and $\varpi:=\left(\omega_1\quad\omega_2\quad \ldots \quad \omega_m\right)^\tr$ be the matrix presentation of this frame, where $\tr$ refers to the matrix transpose. Then we define the \emph{matrix of Gauss-Manin connection}, which is denoted by $\GM_\varpi$, as follow
\[
\nabla \varpi:=\left(\nabla\omega_1\quad\nabla\omega_2\quad \ldots \quad \nabla\omega_m\right)^\tr=\GM_\varpi \otimes \varpi.
\]
We are noting that for any $z\in P$ and any $j\in\{1,2,\ldots,m\}$, $\omega_j(z)\in H^k_\dR(W_z;\C)$ and we can present it by a $k$-form on $W_z$ that we denote it also by $\omega_j(z)$.
}
\end{rem}

Each fiber $H^k_\dR(W_z;\C)$ of
$H^k_\dR(\W/P)$ has a Hodge filtration, and this yields a
decreasing filtration of $H^k_\dR(\W/P)$ by holomorphic
subbundles
\begin{equation}\label{eq filtration of family}
\F^\bullet H^k:\,\,\{0\}={\cal F}^{k+1}\subset {\cal F}^k\subset\ldots\subset {\cal F}^1\subset {\cal F}^0=H^k_\dR (\W/P),
\end{equation}
such that for any $z\in P$ and any $p\in \{0,1,2,\ldots,k\}$
\[
{\cal F}^p_z\cong F^p(W_z)=\bigoplus\limits_{p\leq r\leq k}H^{r,k-r}(W_z).
\]
The filtration $\F^\bullet H^k$ given in (\ref{eq
filtration of family}), is also called \emph{Hodge filtration of}
$H^k_\dR(\W/P)$.
\begin{theo}\label{theo griffiths t} \textrm{{\bf (Griffiths transversality)}}
Under above terminologies, following holds:
\[
\nabla {\cal F}^p \subset \Omega_P^1 \otimes {\cal F}^{p-1},\,\ p=1,2,\ldots k.
\]
\end{theo}

\subsection{Picard-Fuchs Equation}

Here we consider the Hodge filtration $\F^\bullet H^n$ of $H^n(\W/P)$ and fix the
local section $\omega \in \F^n$; indeed for any $z\in P$, $\omega(z)\in H^{n,0}(W_z)$ is a holomorphic $n$-form. Let
$\cal D$ be the ring of linear differential operators on $P$.
If $\dim P=r$ and $z_1,z_2,\ldots,z_r$ is a local coordinate of
$(P,0)$, then we have ${\cal D}=\C(z_1,z_2,\ldots,z_r)[\partial_1,\partial_2,\ldots,\partial_r]$,
where $\C(z_1,z_2,\ldots,z_r)$ is the ring of convergent
power series of $z_1,z_2,\ldots,z_r$ and
$\partial_i=\frac{\partial}{\partial z_i}$. We define the
$\O_P$-homomorphism $\Psi:{\cal D}\to H_\dR^n(\W/P)$, which for
vector fields $v_1,v_2,\ldots,v_k$ on $P$ is determined by
\[
\Psi(v_1v_2\ldots v_k)=\nabla_{v_1}\nabla_{v_2}\ldots \nabla_{v_k} \omega.
\]
By this definition, $\Psi$ gives the structure of a $\cal
D$-module to $H_\dR^n(\W/P)$.
\begin{defi}
{\rm The ideal $\mathcal{I}=\ker \Psi$, consist of differential operators that annihilate $\omega$, by definition
is called \emph{Picard-Fuchs ideal} and any $L\in \mathcal{I}$
is called a \emph{Picard-Fuchs equation}. }
\end{defi}

\begin{ass}\rm
In what follows in this section, we suppose that $\W$ is a one parameter family of
$n$-dimensional compact K\"ahler manifolds, i.e., $\dim P=1$.
\end{ass}

Let $z$ be a
coordinate of $(P,0)$ and define the differential operator
$\vartheta:=\nabla_{z\frac{\partial}{\partial z}}$. Then
$(H^n_\dR(\W/P),\vartheta)$ is a differential $\C(z)$-module.
Considering the terminologies introduced in \S \ref{section
do}, we present the following definition of Picard-Fuchs
equation.

\begin{defi}
{\rm Let $\W$ be a one parameter family of $n$-dimensional
compact K\"ahler manifolds and $\omega\in H^n_\dR(\W/P)$ be a
fixed non-zero element. Then the minimal operator of $\omega$
is called the \emph{Picard-Fuchs equation associated with} $\omega$. }
\end{defi}

\begin{ass} \label{ass 0}\rm
From now on, we suppose that there exists a nowhere vanishing holomorphic $n$-form $\omega\in \F^n$ such that the  Picard-Fuchs equation $\L$  associated with it is of order $n+1$ given as follow
  \begin{equation}\label{eq pf n 2-4}
\L= \vartheta^{n+1}-a_n(z)\vartheta^n-\ldots-a_1(z)\vartheta-a_0(z),
\end{equation}
where $a_i(z)\in \Q(z)$, $i=0,1,\ldots,n$. Therefore, by definition $\L\omega=0$.
\end{ass}

\subsection{Intersection Form}

For any $\alpha,\xi\in H^n_\dR(\W/P)$, the
\emph{intersection form} of $\alpha$ and $\xi$ by definition
is
$$
\langle \alpha,\xi \rangle(z):={\rm Tr}(\alpha(z) \smallsmile \xi(z)),\,\, \forall z\in P,
$$
in which "$\smallsmile$" refers to the cup product. In de
Rham cohomology, the cup product of differential forms is
induced by the wedge product, hence in the family $\W$ the intersection form is defined as follow
\begin{equation}\label{eq intfor}
\langle \alpha,\xi \rangle(z)=\int_{W_z}\alpha(z)\wedge \xi(z).
\end{equation}

We state below a lemma that follows easily from properties of wedge product.
\begin{lemm} \label{lemm if property}
Followings hold:
\begin{description}
  \item[(i)] $\langle \alpha,\xi \rangle=(-1)^n\langle
      \xi,\alpha \rangle$, for any $\alpha,\xi \in
      H_\dR^n(\W/P)$.
  \item[(ii)] If $\F^\bullet H^n$ is the Hodge filtration of
      $H_\dR^n(\W/P)$, then
  \begin{equation}\label{eq gt n}
\langle \F^i,\F^j\rangle=0, \textit{ for } \ i+j\geq n+1.
\end{equation}
\end{description}
\end{lemm}

\subsection{Self-Duality} \label{subsection sd}

Here we give the proof of Proposition \ref{prop 1}. First we fix following notation.

\begin{nota}
{\rm By notation, for $i=1,\ldots,n+1$, we define $\omega_i:=\vartheta^{i-1} \omega$.
}
\end{nota}

We know that $\omega_1=\omega\in \F^n$, hence by Griffiths transversality $\omega_i\in \F^{(n+1)-i}$. Therefore, Lemma \ref{lemm if property}{\bf (ii)} implies that \[
\langle \omega_1,\omega_i \rangle=0,\,\, i=1,2,\ldots,n.
\]
One can find in \cite[\S~4.5]{batst} that
\begin{equation} \label{eq omomn+1}
\langle \omega_1,\omega_{n+1} \rangle(z) = c_0\exp \left( {\frac{2}{n+1} \int_0^z
a_n(v) \frac{dv}{v}}\right),
\end{equation}
for some nonzero constant $c_0$. If we denote by $\tilde{a}(z):= c_0\exp \left( {\frac{2}{n+1}
\int_0^z a_n(v) \frac{dv}{v}}\right)$, then for any $i\in
\{1,\ldots,n\}$
\begin{equation}\label{eq codiagonal}
\langle \omega_i,\omega_{n+2-i} \rangle=(-1)^{i-1}\tilde{a}.
\end{equation}
To see this, first note that by Lemma \ref{lemm if property}{\bf (ii)} we have $\langle \omega_{j+1},\omega_{n-j} \rangle=0,\,\, j=0,1,\ldots,n-1$. On the other hand we know that
$$\vartheta \langle \omega_{j+1},\omega_{n-j} \rangle=\langle \vartheta \omega_{j+1},\omega_{n-j} \rangle+\langle \omega_{j+1},\vartheta \omega_{n-j} \rangle=\langle \omega_{j+2},\omega_{n-j} \rangle+\langle \omega_{j+1},\omega_{n-j+1} \rangle=0,$$
where in the first side of above equation by $\vartheta$ we mean the usual derivation operator $z\frac{\partial}{\partial z}$.
Thus we obtain $\langle \omega_{j+2},\omega_{n-j} \rangle=-\langle \omega_{j+1},\omega_{n-j+1} \rangle$, from which follows (\ref{eq codiagonal}).

\begin{prop} \label{prop h^ij=1}
Let $\F^\bullet H^n$ be the Hodge filtration of $H^n_\dR(\W/P)$. Then $\dim H^n_\dR(\W/P)=n+1$ if and only
if $\dim \F^i/\F^{i+1}=1$ for any $i\in \{0,1,\ldots,n\}$.
\end{prop}
{\bf Proof.} If $\dim \F^i/\F^{i+1}=1$, then it is evident  that
$\dim H^n_\dR(\W/P)=n+1$. Conversely suppose that $\dim
H^n_\dR(\W/P)=n+1$. Then it is enough to prove that $\dim
\F^i/\F^{i+1}\neq 0$. By (\ref{eq omomn+1}) we know that $\langle
\omega_1,\omega_{n+1}\rangle\neq 0$, hence Lemma \ref{lemm if property}{\bf (ii)} implies that $\omega _{n+1}\in \F^0\setminus \F^1$. Now to prove $\dim \F^i/\F^{i+1}\neq 0$, by contradiction suppose that there is a $j\in \{1,2,3,\ldots,n-1\}$ such that $\dim \F^j/\F^{j+1}= 0$, and hence $\F^{j+1}=\F^j$. We know that $\omega_{(n+1)-j}\in \F^j$, thus by Griffiths transversality $\omega_{(n+1)-j+1}=\vartheta \omega_{(n+1)-j}\in \F^{j+1}=\F^j$. Again by using of Griffiths transversality we obtain that $\omega_{(n+1)-j+2}\in \F^j$. By continuing this process it follows that $\omega_{n+1}\in \F^j$, which contradicts $\omega _{n+1}\in \F^0\setminus \F^1$. \hfill\(\blacksquare\)

\begin{ass} \label{ass 1}
{\rm In the rest of this section we assume that for any $i\in \{0,1,\ldots,n\}$, $\dim
      \F^i/\F^{i+1}=1$, or equivalently $\dim H^n_\dR(\W/P)=n+1$.
 }
\end{ass}

\begin{rem}\label{rem nonzero atild}
{\rm
Assumption \ref{ass 1} yields that $\dim \F^i=(n+1)-i,\,\ i=0,1,\ldots,n+1$. It is equivalent to say $\dim \h^{i,j}(W_z)=1$
      for any $z\in P$ and any non-negative integers $i,j$
      with $i+j=n$.
}
\end{rem}

\begin{prop}\label{theo frame}
The set $\{\omega_1,\omega_2,\ldots,\omega_{n+1}\}$ construct a
frame for $H_\dR^n(\W/P)$ such that for any $i \in
\{1,2,\ldots,n+1\}$,
\begin{equation}\label{eq partial 1127}
\omega_i \in \F^{(n+1)-i}\setminus \F^{(n+2)-i}.
\end{equation}
\end{prop}
{\bf Proof. }We know that $\dim H^n_\dR(\W/P)=n+1$, hence it is
enough to show that for any $z$, the set
$\{\omega_1(z),\omega_2(z),\ldots,\omega_{n+1}(z)\}$ is
linearly independent. To this end, suppose that there are
constants $b_1,b_2,\ldots,b_{n+1}$ such that
$
b_1\omega_1(z)+b_2\omega_2(z)+\ldots+b_{n+1}\omega_{n+1}(z)=0.
$
If we set
$
k:={\rm max}\{i\,|\,b_i\neq 0,\, i=1,2,\ldots,n+1\},
$
then we can write
\[
\omega_k(z)=c_1\omega_1(z)+c_2\omega_2(z)+\ldots+c_{k-1}\omega_{k-1}(z),
\]
in which $c_i=\frac{b_i}{b_k}$. By intersecting
$\omega_k$ with $\omega_{n+2-k}$, and using Lemma \ref{lemm if property}{\bf (ii)} we have
\[
\langle \omega_k,\omega_{n+2-k}\rangle(z)=c_1(z)\langle \omega_1,\omega_{n+2-k}\rangle(z)+\ldots+c_{k-1}(z)\langle \omega_{k-1},\omega_{n+2-k}\rangle(z)=0.
\]
On account of (\ref{eq codiagonal}) we get $\langle \omega_k,\omega_{n+2-k}\rangle(z)\neq 0$, which is an contradiction. Thus for any $i\in \{1,2,\ldots,n+1\}$, $b_i=0$.\\
To prove (\ref{eq partial 1127}), first note that Griffiths transversality implies that $\omega_i \in \F^{(n+1)-i}$,
$i=1,2,\ldots,n+1$. On the other hand, since $\dim
\F^{(n+2)-i}=i-1$ and $\{\omega_1,\omega_2,\ldots,\omega_i\}$
is an independent subset of $H^n_\dR(\W/P)$, it follows that $\omega_i \notin
\F^{(n+2)-i}$. \hfill\(\blacksquare\)\\

Finally in the following proposition we give the proof of Proposition \ref{prop 1}.

\begin{prop}\label{theo pp}
The picard-Fuchs equation $\L$ satisfies the property {\rm (P)},
or equivalently $\L$ is self-dual.
\end{prop}
{\bf Proof. } Consider the intersection form defined as follow
\[
\langle .,. \rangle: H_\dR^n(\W/P)\times H_\dR^n(\W/P)\to \C(z).
\]
Equation (\ref{eq codiagonal}) implies that $\langle .,.
\rangle$ is non-degenerate, and the Lemma \ref{lemm if
property}{\bf (i)} verifies that $\langle .,. \rangle$ is a
$(-1)^n$-symmetric form. Lemma \ref{lemm if property}{\bf
(ii)} guaranties that in frame
$\{\omega,\vartheta\omega,\ldots,\vartheta^n\omega\}$ we have,
\[
\langle \omega,\vartheta^i\omega \rangle=0,\,\, {\rm for }\,\, i=0,1,\ldots,n-1.
\]
Hence by Definition \ref{defi property p}, $\L$ satisfies the property (P). This is equivalent with self-duality of $\L$ by Proposition \ref{prop bog}.\hfill\(\blacksquare\)\\

\section{Darboux-Halphen-Ramanujan Vector Field} \label{section dhr}

In this section  $\pi:\W \to P$ refers to a one parameter family of $n$-dimensional \cy manifolds, or equivalently it is a complex deformation of an $n$-dimensional \cy manifold $W:=W_0$. Let $z$ be a local coordinate of $(P,0)$. Then for any $z\in P$, $W_z$ is a \cy $n$-fold. We know that Calabi-Yau manifold $W$ is a compact K\"{a}hler manifold whose, up to multiplication by a
constant, has a unique nowhere vanishing holomorphic $n$-form $\omega\in H^{n,0}(W)$. Thus, there is a holomorphic section of $H^n_\dR(\W/P)$ that at $0$ coincides with $\omega$ and we denote it also by $\omega$. Hence $\omega\in \F^n$, where $\F^\bullet H^n$ is the Hodge filtration of $H^n_\dR(\W/P)$, and $\omega(z)\in H^{n,0}(W_z)$ is a nowhere vanishing holomorphic $n$-form of $W_z$ for any $z$. In this section we fix $\omega\in \F^n$ and suppose that $\omega$ satisfies Assumption \ref{ass 0} and $\F^\bullet H^n$ satisfies Assumption \ref{ass 1}. Throughout this section, we employ the notations of pervious sections.
\begin{exam}\label{exam 14 cases}
By now, as I know, there are $14$ examples of one parameter families of \cy 3-folds satisfying the hypothesis of $\W$ given above.
Any of these $14$ families is mirror symmetry of $14$ structures given in {\rm Table \ref{table 14 case}}. In this table $X(d_1,d_2,\ldots,d_r)\subset \P^s(l_0,l_1,\ldots,l_s)$ refers to the complete intersection of $r$ hypersurfaces of degree $d_1,d_2,\ldots,d_r$ in weighted projective space $\P^s(l_0,l_1,\ldots,l_s)$ with $r\leq s$, such that $\sum _{i=1}^r d_i=\sum _{j=0}^s l_j$. The Picard-Fuchs equation $\L$ associated with the nowhere vanishing holomorphic $3$-form of any of these families is a hypergeometric equation given as follow:
\begin{equation}\label{eq pf3 14}
\L=\vartheta^4-cz(\vartheta+r_1)(\vartheta+r_2)(\vartheta+1-r_2)(\vartheta+1-r_1),
\end{equation}
where $r_1,r_2,c$ are given in {\rm Table \ref{table 14 case}}. Note that $\sharp\, 1$ is the family of quintic 3-folds that we pointed it out in \S \ref{section int}.  For more details one can see the
references given in {\rm Table \ref{table 14 case}}.
\begin{table}[hbt!]\tiny
\centering
  \begin{tabular}{|c||c|c|c||l||c|}
    \hline
    $\sharp$ & $r_1$ & $r_2$ & $c$ & Structure & References \\
    \hline\hline
    $1$ & $1/5$ & $2/5$ & $5^5$ & $X(5)\subset \P^4$ & \cite{can91,greple} \\
    \hline
    $2$ & $1/6$ & $2/6$ & $2^53^6$ & $X(6)\subset\P^4(2,1,1,1,1)$ & \cite{mo92} \\
    \hline
    $3$ & $1/8$ & $3/8$ & $2^{18}$ & $X(8)\subset\P^4(4,1,1,1,1)$ & \cite{mo92} \\
    \hline
    $4$ & $1/10$ & $3/10$ & $2^95^6$ & $X(10)\subset\P^4(5,2,1,1,1)$ & \cite{mo92} \\
    \hline
    $5$ & $1/3$ & $1/3$ & $3^{6}$ & $X(3,3)\subset\P^5$ & \cite{libtei} \\
    \hline
    $6$ & $1/4$ & $2/4$ & $2^{10}$ & $X(2,4)\subset\P^5$ & \cite{libtei} \\
    \hline
    $7$ & $1/3$ & $1/2$ & $2^43^3$ & $X(2,2,3)\subset\P^6$ & \cite{libtei} \\
    \hline
    $8$ & $1/2$ & $1/2$ & $2^8$ & $X(2,2,2,2)\subset\P^7$ & \cite{libtei} \\
    \hline
    $9$ & $1/4$ & $1/4$ & $2^{12}$ & $X(4,4)\subset\P^5(2,2,1,1,1,1)$ & \cite{klethe} \\
    \hline
    $10$ & $1/6$ & $1/6$ & $2^83^6$ & $X(6,6)\subset\P^5(3,3,2,2,1,1)$ & \cite{klethe} \\
    \hline
    $11$ & $1/4$ & $1/3$ & $2^63^3$ & $X(3,4)\subset\P^5(2,1,1,1,1,1)$ & \cite{klethe} \\
    \hline
    $12$ & $1/6$ & $3/6$ & $2^83^3$ & $X(2,6)\subset\P^5(3,1,1,1,1,1)$ & \cite{klethe} \\
    \hline
    $13$ & $1/6$ & $1/4$ & $2^{10}3^3$ & $X(4,6)\subset\P^5(3,2,2,1,1,1)$ & \cite{klethe} \\
    \hline
    $14$ & $1/12$ & $5/12$ & $12^6$ & $X(2,12)\subset\P^5(6,4,1,1,1,1)$ & \cite{dormor} \\
    \hline
  \end{tabular}
 \caption{\tiny \cy 3-folds \cite{Chen08monodromyof}}
  \label{table 14 case}
\end{table}
\end{exam}

Hodge filtration of $H^n_{\dR}(\W/P)$ is as follow
\begin{equation}\label{hodgefilteration n}
\F^\bullet H^n:\,\, \{0\}=\F^{n+1}\subset \F^n\subset \ldots \subset \F^1\subset \F^0=H^{n}_{\dR} (\W/P),\ \ \dim \F^i=(n+1)-i,
\end{equation}
and as we saw in Proposition \ref{theo frame}, $\{\omega_1,\omega_2,\ldots,\omega_{n+1}\}$ construct a frame of $H^n_{\dR}(\W/P)$ such that
\begin{equation}\label{eq 3714 1}
\omega_i \in \F^{(n+1)-i}\setminus \F^{(n+2)-i}.
\end{equation}
By using of Picard-Fuchs equation (\ref{eq pf n 2-4}) we have
\begin{equation}\label{eq vtheta om}
\vartheta^{n+1}\omega=\vartheta\omega_{n+1}=a_0\omega_1+a_1\omega_2+\ldots+a_n\omega_{n+1}.
\end{equation}
Hence, considering Remark \ref{obs gm mat}, if we apply the Gauss-Manin connection to the column of $n$-forms
$
\varpi=\left( {\begin{array}{*{20}c}
   {\omega _1 } & {\omega _2 } &  \ldots  & {\omega _{n + 1} }  \\
\end{array}} \right)^\tr,
$
then
\begin{equation}\label{eq gmomega n}
\nabla \varpi=\GM_\varpi\otimes\varpi,
\end{equation}
where
\begin{equation}\label{eq gm varpi n}\small
\GM_\varpi=\frac{1}{z}\left(
         \begin{array}{cccccc}
           0 & 1 & 0 & \ldots & 0 & 0 \\
           0 & 0 & 1 & \ldots & 0 & 0 \\
           \vdots & \vdots & \vdots & \ddots & \vdots & \vdots \\
           0 & 0 & 0 & \ldots & 1 & 0 \\
           0 & 0 & 0 & \ldots & 0 & 1 \\
           a_0 & a_1 & a_2 & \ldots & a_{n-1} & a_n \\
         \end{array}
       \right)dz.
\end{equation}
To see this, for $j=1,2,\ldots,n$, we have
\[
z\nabla_{\frac{\partial}{\partial z}}\omega_j=\nabla_{z\frac{\partial}{\partial z}}\omega_j=\vartheta \omega_j=\omega_{j+1}\Longrightarrow \, \nabla \omega_j=\frac{1}{z}dz\otimes \omega_{j+1}.
\]
Analogously, on account of (\ref{eq vtheta om}), for $\omega_{n+1}$ we obtain
$
\nabla \omega_{n+1}=\frac{1}{z}\sum\limits_{i=0}^{n}a_i dz\otimes \omega_{i+1}.
$

In this section we are going to prove Theorem \ref{theo 1 int}. In order to do this, we will treat with intersection form, but because of different behaviors of intersection form for odd or even integer $n$,  see Lemma \ref{lemm if property}{\bf (i)}, we separate the cases for odd and even
integers. First, we state the results in the odd case in \S \ref{section cy n}. In particular for $n=3,5$ we give an explicit computation of results in \S\ref{section 3-dim} and \S\ref{section 5-dim}.

\subsection{Odd Case} \label{section cy n}
In the whole of this subsection $n$ is considered to be an odd
positive integer. If we define the intersection form matrix as follow
\begin{equation} \label{eq int matrix n}
\Omega=\left(\Omega_{ij}\right)_{1\leq i,j\leq n+1}:=\left( \langle\omega_i,\omega_j\rangle \right)_{1\leq i,j\leq n+1},
\end{equation}
then Lemma \ref{lemm if property}{\bf (i)} implies that $\Omega^\tr=-\Omega$, and hence $\Omega_{ii}=0, \, i=1,2,\ldots,n+1$. Lemma \ref{lemm if property}{\bf (ii)} yields $\Omega_{ij}=\langle\omega_i,\omega_j\rangle=0$ for $i+j\leq n+1$, and by (\ref{eq codiagonal}) we find $\Omega_{i({n+2-i})}=(-1)^{i-1}\tilde{a}$ for any $i=1,2,\ldots,n+1$. Therefore, we can state the matrix $\Omega$ as follow
\begin{equation}\label{eq omega n}\small
\Omega=\left( {\begin{array}{*{20}c}
   0 & 0 &  \ldots   & 0 & {\tilde a}  \\
   0 & 0 &  \ldots   & { - \tilde a} & {\Omega _{2(n + 1)} }  \\
       \vdots  &  \vdots  &  {\mathinner{\mkern2mu\raise1pt\hbox{.}\mkern2mu
 \raise4pt\hbox{.}\mkern2mu\raise7pt\hbox{.}\mkern1mu}}  &  \vdots  &  \vdots   \\
   0 & {\tilde a} &  \ldots  & 0 & {\Omega _{n(n + 1)} }  \\
   { - \tilde a} & {-\Omega _{2(n + 1)} } &  \ldots  & {-\Omega _{n(n + 1)} } & 0  \\
\end{array}} \right).
\end{equation}

\begin{defi}\rm
We say that a basis
$\{\alpha_1,\alpha_2,\ldots,\alpha_{n+1}\}$ of
$H^n_{\dR}(W;\C)$ is \emph{compatible with its Hodge filtration}, if for any $i \in
\{1,2,\ldots,n+1\}$
\begin{equation} \label{eq 3714 2}
\alpha_i\in F^{n+1-i}\setminus F^{n+2-i}.
\end{equation}
\end{defi}

Next we introduce a special moduli space of \cy
manifold $W$ that in the rest of this text will be in interest. To do this, we first
provide an equivalence  relation.

\begin{defi} \rm
Let $W_1,W_2$ be two \cy $n$-folds and
$\{\alpha_1^i,\alpha_2^i,\ldots,\alpha_{n+1}^i\}$ be a basis of
$H_\dR^n(W_i;\C), i=1,2$, compatible with its Hodge filtration.
Then we write
\begin{equation} \label{eq equi rel}
(W_1,[\alpha_1^1,\alpha_2^1,\ldots,\alpha_{n+1}^1])\sim (W_2,[\alpha_1^2,\alpha_2^2,\ldots,\alpha_{n+1}^2])
\end{equation}
if and only if there exist a biholomorphism $\varphi:W_1\to
W_2$ such that $\varphi^\ast(\alpha_j^2)=\alpha_j^1,\,\
j=1,2,\ldots,n+1$. It is obvious that "$\sim$" is an
equivalence relation. For the \cy $n$-fold $W$, and a basis
$\{\alpha_i\}_{i=1}^{n+1}$ of $H_\dR^n(W;\C)$ compatible with
its Hodge filtration, the \emph{moduli space} $\tilde{\Tm}$ of
pair $(W,[\alpha_1,\alpha_2,\ldots,\alpha_{n+1}])$ is defined
under above equivalence relation (\ref{eq equi rel}).
\end{defi}

\begin{rem} \rm
We know that the family $\pi:\W \to P$ is the complex
deformation of $W$. Hence for any different $z_1, z_2\in P$,
$W_{z_1}$ and $W_{z_2}$ are not biholomorph. We thus have two different members
$(W_{z_1},[\omega_1(z_1),\omega_2(z_1),\ldots\omega_{n+1}(z_1)])$
and
$(W_{z_2},[\omega_1(z_2),\omega_2(z_2),\ldots\omega_{n+1}(z_2)])$
of moduli space $\tilde{\Tm}$. Also suppose
that $\{\mu_i\}_{i=1}^{n+1}$ and $\{\nu_i\}_{i=1}^{n+1}$ are
two bases of $H_\dR^n(W;\C)$ compatible with
its Hodge filtration. If for
any
\[
\varphi\in \textrm{Aut}(W):=\{f:W\to W|\, f \textrm{ is a biholomorphism}\},
\]
it does not preserve the bases, i.e., there exist a
$j\in \{1,2,\ldots,n+1\}$ such that $\varphi^\ast \nu_j\neq
\mu_j$, then $(W,[\mu_1,\mu_2,\ldots,\mu_{n+1}])$ and
$(W,[\nu_1,\nu_2,\ldots,\nu_{n+1}])$ yield two different
elements of moduli space $\tilde{\Tm}$.
\end{rem}

As we fixed in the beginning of this section, for any $z\in P$,
$\{\omega_1(z),\omega_2(z),\ldots,\omega_{n+1}(z)\}$ construct
a basis for $H^n_\dR(W_z,\C)$ that is compatible with its Hodge
filtration. By abuse of notation, we remove the letter $z$
from this basis and denote it by
$\{\omega_1,\omega_2,\ldots,\omega_{n+1}\}$, and hence
$(W_z,[\omega_1,\omega_2,\ldots\omega_{n+1}])\in \tilde{\Tm}$.
Let $S$ be the change of basis matrix $\alpha=S\varpi$, where $\{\alpha_i\}_{i=1}^{n+1}$ is a basis of $H_\dR^n(W_z;\C)$ compatible with its Hodge
filtration, and
$
\alpha=\left( {\begin{array}{*{20}c}
   {\alpha _1 } & {\alpha _2 } &  \ldots  & {\alpha _{n + 1} }  \\
\end{array}} \right)^\tr.
$
Then (\ref{eq 3714 1}) and (\ref{eq 3714 2}) imply that $S$ is a lower triangular matrix which we consider it as follow
\begin{equation}\label{eq s n}\small
S=\left(
         \begin{array}{ccccc}
           s_{11} & 0 & 0 & \ldots & 0 \\
           s_{21} & s_{22}& 0 & \ldots & 0 \\
           s_{31} & s_{32} & s_{33} & \ldots & 0 \\
           \vdots & \vdots & \vdots & \ddots & \vdots  \\
           s_{(n+1)1} & s_{(n+1)2} & s_{(n+1)3} & \ldots & s_{(n+1)(n+1)}  \\
              \end{array}
       \right).
\end{equation}
Hence the entries of $S$ present
coordinates of a chart of $\tilde{\Tm}$ that we will employ it soon.

\begin{lemm} \label{lemm phi gm alpha}
Let $\{\alpha_i\}_{i=1}^{n+1}$ be a frame of
$H^n_{\dR}(\W/P)$ compatible with its Hodge filtration.
\begin{description}
  \item[(i)] If we define $\Psi:=\left(
      \langle\alpha_i,\alpha_j\rangle \right)_{1\leq
      i,j\leq n+1}$, then
$
\Psi=S\Omega S^\tr.
$
  \item[(ii)] If we set
      $\nabla\alpha=\GM_{\alpha}\otimes\alpha$, then
\begin{equation}\label{eq gm alpha}
\GM_\alpha=(dS + S.\GM_\varpi){S^{ - 1}},
\end{equation}
where
\begin{equation}\label{eq s}\small
dS=\left(
         \begin{array}{ccccc}
           ds_{11} & 0 & 0 & \ldots & 0 \\
           ds_{21} & ds_{22}& 0 & \ldots & 0 \\
           ds_{31} & ds_{32} & ds_{33} & \ldots & 0 \\
           \vdots & \vdots & \vdots & \ddots & \vdots  \\
           ds_{(n+1)1} & ds_{(n+1)2} & ds_{(n+1)3} & \ldots & ds_{(n+1)(n+1)}  \\
              \end{array}
       \right).
\end{equation}
\end{description}
\end{lemm}
{\bf Proof.}\begin{description}
              \item[(i)] By using of $\alpha=S\varpi$, verifying $\Psi=S\Omega S^\tr$ is an easy exercise of linear algebra.
              \item[(ii)] If we apply the Gauss-Manin connection to the equation  $\alpha  = S\varpi$, and considering $\nabla\varpi=\GM_{\varpi}\otimes\varpi$, then we have
\[\begin{array}{l}
  \nabla \alpha  = dS \otimes \varpi  + S\nabla \varpi  = (dS + S.\GM_\varpi) \otimes \varpi \\
\,\,\,\,\quad\, = (dS + S.\GM_\varpi){S^{ - 1}} \otimes \alpha,
\end{array}\]
which completes the proof. \hfill\(\blacksquare\)
            \end{description}

Following proposition give a more important step of the proof
of Theorem \ref{theo 1 int}.

\begin{prop}\label{prop h}
Let {\rm $\tilde{\textsf{T}}$} be the moduli of
$(W,[\alpha_1,\alpha_2,\ldots,\alpha_{n+1}])$, where $\{\alpha_i\}_{i=1}^{n+1}$ is a basis of
$H^n_{\dR}(W;\C)$ compatible with its Hodge filtration.
Then there exist a unique vector field $\tilde{\H}$ and unique
meromorphic functions $y_i$, $i=1,2,\ldots,n-2$, on {\rm
$\tilde{\Tm}$} such that
\begin{equation}\label{eq gm halphen}
\nabla_{\tilde{\H}}\alpha=Y\alpha,
\end{equation}
in which $\alpha=\left( {\begin{array}{*{20}c}
   {\alpha _1 } & {\alpha _2 } &  \ldots  & {\alpha _{n + 1} }  \\
\end{array}} \right)^\tr$, and $Y$ is given by {\rm (\ref{eq y n})}.
\end{prop}

{\bf Proof.} The idea of the proof is to present the vector
field $\tilde{\H}$ explicitly in a chart of $\tilde{\Tm}$. It
is easily seen that the dimension of $\tilde{\Tm}$ is $k+1$,
where $k=\frac{(n+1)(n+2)}{2}$.  For any
$(W_z,[\alpha_1,\alpha_2,\ldots,\alpha_{n+1}])\in \tilde{\Tm}$,
let $S$ be the change of basis matrix $\alpha=S\varpi$ given in
(\ref{eq s n}). We consider the chart
$t=(t_0,t_1,\ldots,t_{k})$ of $\tilde{\Tm}$, for which the
coordinates are defined as $t_0=z, t_1=s_{11},
t_2=s_{12},\ldots, t_{k}=s_{(n+1)(n+1)}$. We suppose that the vector field
$\tilde{\H}$ is given as follow
\[
\tilde{\H}=\sum_{i=0}^{k}\tilde{\H}_i(t)\frac{\partial}{\partial t_i},
\]
where $\tilde{\H}_i$'s, $i=0,1,\ldots,k$, are meromorphic
functions on $\tilde{\Tm}$. Since $\tilde{\H}$ satisfies
$\nabla_{\tilde{\H}}\alpha=Y\alpha$, Lemma \ref{lemm phi
gm alpha}{\bf (ii)} implies that
\begin{equation}\label{eq 3914 1}
(dS + S.\GM_\varpi){S^{ - 1}}({\tilde{\H}})=Y.
\end{equation}
We have  $S.\GM_\varpi(\tilde{\H})=\dot{z}S.\widehat{\GM}_\varpi$, where $\dot{z}(t):=\tilde{\H}_0(t)$ and $\widehat{\GM}_\varpi$ is defined by $\GM_\varpi=\widehat{\GM}_\varpi dz$.
Also if we define $\dot{s}_{11}(t):=\tilde{\H}_1(t),\, \dot{s}_{21}(t):=\tilde{\H}_2(t),\, \ldots, \, \dot{s}_{(n+1)(n+1)}(t):=\tilde{\H}_k(t)$,
then we have $dS(\tilde{\H})=\dot{S}$, where
\begin{equation}\label{eq dots n}
\dot{S}=\left(
         \begin{array}{ccccc}
           \dot{s}_{11} & 0 & 0 & \ldots & 0 \\
           \dot{s}_{21} & \dot{s}_{22}& 0 & \ldots & 0 \\
           \dot{s}_{31} & \dot{s}_{32} & \dot{s}_{33} & \ldots & 0 \\
           \vdots & \vdots & \vdots & \ddots & \vdots  \\
           \dot{s}_{(n+1)1} & \dot{s}_{(n+1)2} & \dot{s}_{(n+1)3} & \ldots & \dot{s}_{(n+1)(n+1)}  \\
              \end{array}
       \right).
\end{equation}
Therefore (\ref{eq 3914 1}) gives
$
(\dot{S}+\dot{z}S.\widehat{\GM}_\varpi){S^{ - 1}}=Y,
$
which yields
\begin{equation}\label{eq dots}
\dot{S}=YS-\dot{z}S.\widehat{\GM}_\varpi.
\end{equation}
Consequently we can find $\dot{z}$ (or $\tilde{\H}_0$) and $y_i$'s (that we state them
in Lemma \ref{lemm y} below). Hence all the terms of the right hand side of  (\ref{eq dots}) are determined, from which we can
find $\tilde{\H}_i$'s, $i=1,2,\ldots,k$. Thus, the existence of vector field $\tilde{\H}$ that satisfies (\ref{eq gm halphen}) is verified. The uniqueness of $\tilde{\H}$ and $y_i$'s follow from Lemma \ref{lemm y}{\bf (i)},{\bf (ii)}. \hfill\(\blacksquare\)\\

The proof of Proposition \ref{prop h}, implies more results about entries of $Y$ that we express them in a lemma below. Before that, we provide the following fact as a remark.

\begin{rem} \label{rem siineq0}\rm
Since the matrix $S$ is the change of basis matrix, it is invertible; thus for any $1\leq i\leq n+1$, $s_{ii}\neq 0$.
\end{rem}

\begin{lemm}\label{lemm y}
The equation \begin{equation}\label{eq dots2}
\dot{S}=YS-\dot{z}S.\widehat{\GM}_\varpi,
\end{equation}
 implies that:
\begin{description}
  \item[(i)] $\dot{z}=\frac{zs_{22}}{s_{11}}=\frac{zs_{(n+1)(n+1)}}{s_{nn}}$.
  \item[(ii)]
      $y_{i-1}=\frac{s_{22}s_{ii}}{s_{11}s_{(i+1)(i+1)}}$,
      for all $i=2,3,\ldots, n-1$.
  \item[(iii)] Moreover, if $S\Omega S^\tr=\Phi$, then
      $y_{i-1}=-y_{n-i}$, for $i\neq \frac{n+1}{2}$; and
      $$y_{\frac{n-1}{2}}=(-1)^{\frac{n+3}{2}}\frac{\tilde{a}s_{22}s_{\frac{n+1}{2}\frac{n+1}{2}}^2}{s_{11}}.$$
      In the other word
\begin{equation}\label{eq yphi}
Y\Phi=-\Phi Y^\tr.
\end{equation}
\end{description}

\end{lemm}

{\bf Proof.} Let's define
$
B=(b_{ij})_{1\leq i,j\leq n+1}:=YS-\dot{z}S.\widehat{\GM}_\varpi.
$
\begin{description}
  \item[(i)] The equation (\ref{eq dots2}) implies that
  $
  b_{12}=s_{22}-\frac{\dot{z}}{z}s_{11}=0$ and $b_{n(n+1)}=s_{(n+1)(n+1)}-\frac{\dot{z}}{z}s_{nn}=0,
  $
  which prove {\bf (i)}.
  \item[(ii)] The proof of {\bf (ii)} follows from {\bf (i)} and
  $
  b_{i(i+1)}=y_{i-1}s_{(i+1)(i+1)}-\frac{\dot{z}}{z}s_{ii}=0,\,\, i=2,3,\ldots, n-1.
  $
  \item[(iii)] Let's define
$
C=(c_{ij})_{1\leq i,j\leq n+1}:=S\Omega S^\tr.
$
Then equation $C=\Phi$ yields
$
c_{i(n+2-i)}=(-1)^{i+1}\tilde{a}s_{ii}s_{(n+2-i)(n+2-i)}=1,\,\, i=1,2,\ldots,\frac{n+1}{2},
$
from which we obtain
\begin{equation}\label{ eq sii}
s_{(n+2-i)(n+2-i)}=(-1)^{i+1}\frac{1}{\tilde{a}s_{ii}},\,\, i=1,2,\ldots,\frac{n+1}{2}.
\end{equation}
Thus,
$$
\frac{s_{ii}}{s_{(i+1)(i+1)}}=-\frac{s_{(n+1-i)(n+1-i)}}{s_{(n+2-i)(n+2-i)}},\,\, i=1,2,\ldots,\frac{n-1}{2}.
$$
Therefore, on account of {\bf (ii)} the proof of {\bf (iii)} is complete.
\end{description}
 \hfill\(\blacksquare\)

\begin{lemm}\label{lemm varthetaomega}
Let $A:=z\widehat{\GM}_\varpi$. Then following equation holds:
\begin{equation}\label{eq thetaomega}
\vartheta\Omega=A\Omega+\Omega A^\tr.
\end{equation}
\end{lemm}

{\bf Proof.} By using of the fact $\vartheta \langle \omega_i,\omega_j \rangle= \langle \vartheta\omega_i,\omega_j \rangle+\langle \omega_i,\vartheta\omega_j \rangle$ and Picard-Fuchs equation (\ref{eq vtheta om}), the proof is an easy exercise of linear algebra. \hfill\(\blacksquare\)\\

Finally, we are now in a position that can prove Theorem \ref{theo 1 int}.\\

{\bf Proof of Theorem \ref{theo 1 int}.} Let $\tilde{\Tm}$ be the moduli space introduced
in Proposition \ref{prop h}, and suppose that
$(W_z,[\alpha_1,\alpha_2,\ldots,\alpha_{n+1}])\in \tilde{\Tm}$
is an arbitrary element. As we saw in the proof of Proposition \ref{prop
h}, there exist the matrix $S$ such that
$$
\left( \langle\alpha_i,\alpha_j\rangle \right)_{1\leq i,j\leq n+1}=S\Omega S^\tr.
$$
Define the vector subspace $M\subset\textrm{ Mat}_{n+1}(\C)$ to be
$$
M:=\{B=(b_{ij})_{1\leq i,j\leq n+1}\in \textrm{Mat}_{n+1}(\C)|\,\, b_{ij}=0,\,\, if\,\, i\leq n+1-j\}.
$$
If we define the map $f$ as follow
\begin{eqnarray*}
\!\!\!&f:&\!\!\!\tilde{\Tm}\to M\\
\!\!\!&f(&\!\!\!W_z,[\alpha_1,\alpha_2,\ldots,\alpha_{n+1}])=S\Omega S^\tr,
\end{eqnarray*}
then $\Tm=f^{-1}(\Phi)$. Hence \Tm is a subspace of $\tilde{\Tm}$ and to prove the
existence of vector field $\H$ on $\Tm$, it is enough to show
that the vector field $\tilde{\H}$, which was introduced in
Proposition \ref{prop h}, is tangent to \Tm and define
$\H:=\tilde{\H}|_\Tm$. To demonstrate the tangency of
$\tilde{\H}$ to \Tm, it suffices to prove that
$df|_\Tm({\tilde{\H}})=0$, or equivalently verify that
\begin{equation}\label{eq sos}
(\dot{S}\Omega S^\tr+S\dot{\Omega} S^\tr+S\Omega \dot{S}^\tr)|_\Tm=0,
\end{equation}
in which $\dot{\Omega}=d\Omega(\tilde{\H})$. Since $\Omega$ just depends to $z$, it follows that $\dot{\Omega}=\dot{z}\frac{\partial}{\partial z}\Omega$. By using of Lemma \ref{lemm varthetaomega} it is deduced that
 \[
 \dot{\Omega}=\frac{\dot{z}}{z}(A\Omega+\Omega A^\tr)=\dot{z}(\widehat{\GM}_\varpi\Omega+\Omega \widehat{\GM}_\varpi^\tr).
 \]
On the other hand as we saw in (\ref{eq dots}), $\dot{S}=YS-\dot{z}S.\widehat{\GM}_\varpi,$ hence
\[
\dot{S}\Omega S^\tr+S\dot{\Omega} S^\tr+S\Omega \dot{S}^\tr=YS\Omega S^\tr+S\Omega S^\tr Y^\tr.
\]
Since $S\Omega S^\tr|_\Tm=\Phi$, by using of Lemma \ref{lemm
y}{\bf (iii)} we get
\[
(\dot{S}\Omega S^\tr+S\dot{\Omega} S^\tr+S\Omega \dot{S}^\tr)|_\Tm=(YS\Omega S^\tr+S\Omega S^\tr Y^\tr)|_\Tm=(Y\Phi+\Phi Y^\tr)|_\Tm=0,
\]
and the proof of existence of $\H$ is complete.

To prove the uniqueness, first notice that Lemma \ref{lemm
y}(ii) guaranties the uniqueness of $y_i$'s. Hence we just need to
prove that the vector field $\H$ is unique. Suppose that there are
two vector fields $\H_1$ and $\H_2$ such that
$\nabla_{\H_i}\alpha=Y\alpha, \, i=1,2$. If we set ${\rm
R}:=\H_1-\H_2$, then
\begin{equation} \label{eq 4114 1}
\nabla_{{\rm R}}\alpha=0.
\end{equation}
We need to prove that ${\rm R}=0$, and to do this it is enough
to verify that any integral curve of ${\rm R}$ is a constant
point. Assume that $\gamma$ is an integral curve of ${\rm
R}$ given as follow
\begin{align}
\gamma:(\C,&\,0) \to \Tm;\qquad
           x\mapsto \gamma(x). \nonumber
\end{align}
Let's denote $\mathcal{C}:=\gamma(\C,0)\subset \Tm$. We know that the members of
\Tm are in the form of the pairs
$(\widehat{W},[\alpha_1,\alpha_2,\ldots,\alpha_{n+1}])$, where $\widehat{W}$ is a
\cy manifold of the family $\W$, and $\{\alpha_i\}_{i=1}^{n+1}$ form a basis of $H_\dR^n(\widehat{W};\C)$ that is compatible with its Hodge filtration and has constant intersection form matrix $\Phi$. Thus,  for any $x\in (\C,0)$, we
have
$\gamma(x)=(\widehat{W}(x),[\alpha_1(x),\alpha_2(x),\ldots,\alpha_{n+1}(x)])$,
and  the vector field $\rm R$ on $\mathcal{C}$ is reduced to
$\frac{\partial}{\partial x}$ as well. We know that $\widehat{W}(x)$ depends only
on the parameter $z$, and hence $x$ holomorphically depends to $z$.
From this we obtain a holomorphic function $f$ such that
$x=f(z)$. We now proceed to prove that $f$ is constant. Otherwise, by
contradiction suppose that $f'\neq 0$. Then we get
\begin{equation}\label{eq 72414}
\nabla_{\frac{\partial}{\partial x}}\alpha_1= \frac{\partial z}{\partial x}\nabla_{\frac{\partial}{\partial z}}\alpha_1.
\end{equation}
Equation (\ref{eq 4114 1}) gives that
$\nabla_{\frac{\partial}{\partial x}}\alpha_1=0$, but since $\alpha_1=s_{11}\omega_1$, it follows that the right
hand side of (\ref{eq 72414}) is not zero, which
is a contradiction. Thus $f$ is constant and $\widehat{W}(x)$ does not depend
on the parameter $x$.
Since  $\widehat{W}(x)=\widehat{W}$ does not depend on $x$, we can write the Taylor
series of $\alpha_i(x),\,\ i=1,2,3,\ldots, n+1,$ in $x$ at some
point $x_0$ as
$\alpha_i(x)=\sum_j (x-x_0)^j\alpha_{i,j},$
where $\alpha_{i,j}$'s are elements in $H_\dR^n(\widehat{W};\C)$ independent
of $x$. In this way the action of
$\nabla_{\frac{\partial}{\partial x}}$ on $\alpha_i$ is just
the usual derivation $\frac{\partial}{\partial x}$. Again according to
(\ref{eq 4114 1}) we yield $\nabla_{\frac{\partial}{\partial
x}}\alpha_i=0$, and we conclude that $\alpha_i$'s also do not
depend on $x$. Therefore, the image of $\gamma$ is a point.

To prove that
$
\dim \Tm=\frac{(n+1)(n+3)}{4}+1,
$
it is enough to observe that $S\Omega S^\tr=\Phi$ gives
$\frac{(n+1)(n+3)}{4}$ independent equations and that $\W$ is a
one parameter family.
 \hfill\(\blacksquare\)\\

\begin{rem}\label{rem chart t} \rm
Let $m:=\frac{(n+1)(n+3)}{4}$ and fix $m$ entries of $S$. By notation we denote them by $t_1,t_2,\ldots,t_m$ and call them \emph{independent entries} of $S$.
The matrix equation $S\Omega S^\tr=\Phi$ yields $\frac{(n+1)^2}{4}$ independent equations that express the rest of entries of $S$, what we shall call \emph{dependent entries}, in terms of $t_i$'s. For instance, suppose that $s_{ab}$ is a dependent entry of $S$ and $S\Omega S^\tr=\Phi$ gives the expression $s_{ab}=\varphi(t_1,t_2,\ldots,t_m)$. We can obtain $\dot s_{ab}$ in the following two ways:
\begin{description}
  \item[(i)] On account of $s_{ab}=\varphi(t_1,t_2,\ldots,t_m)$, we first get
  \begin{equation} \label{eq sdott}
  \dot s_{ab}=\sum_{i=1}^m \dot t_i\frac{\partial \varphi}{\partial t_i},
  \end{equation}
  and then substitute $\dot t_i$'s from $\dot{S}=YS-\dot{z}S.\widehat{\GM}_\varpi$ in (\ref{eq sdott}).
  \item[(ii)] We first find $\dot s_{ab}$ directly from $\dot{S}=YS-\dot{z}S.\widehat{\GM}_\varpi$, and then using $S\Omega S^\tr=\Phi$ to express $\dot s_{ab}$ just in terms of $t_i$'s.
\end{description}
We say that \emph{the equations $S\Omega S^\tr=\Phi$ and $\dot{S}=YS-\dot{z}S.\widehat{\GM}_\varpi$ are compatible} if {\bf (i)} and {\bf (ii)} give the same result for $\dot s_{ab}$. We are now in a position to introduce a chart of $\Tm$, where $t_i$'s are its coordinates. In order to this, let $t_0:=z$. Then $t=(t_0,t_1,\ldots,t_m)$ gives a chart for $\Tm$ that we will work explicitly with it in \S\ref{section 5-dim} and \S\ref{section 3-dim}.
\end{rem}

The corollary stated below, is an immediate result of Theorem \ref{theo 1 int}.

\begin{coro}
The equations $S\Omega S^\tr=\Phi$ and $\dot{S}=YS-\dot{z}S.\widehat{\GM}_\varpi$ are compatible on $\Tm$.
\end{coro}

Conversely, one can find that the compatibility of equations $S\Omega S^\tr=\Phi$ and $\dot{S}=YS-\dot{z}S.\widehat{\GM}_\varpi$ implies the existence and uniqueness of \DHR vector field. We see this clearly in  \S\ref{section 5-dim} and \S\ref{section 3-dim}, where we compute \DHR vector field explicitly.

\subsection{Even Case} \label{section cy n even}
During this subsection $n$ refers to an even positive integer. As
we mentioned before, the difference of even case with the odd
case is just the symmetry of intersection form. Lemma \ref{lemm
if property} implies that in the odd case the intersection form matrix is anti-symmetric, but in the even case it
is symmetric. Hence, in this section we follow all the notations
and definitions of \S \ref{section cy n}, except the concepts related
with intersection form. In particular the matrix
$\Omega=\left(\Omega_{ij}\right)_{1\leq i,j\leq n+1}:=\left(
\langle\omega_i,\omega_j\rangle \right)_{1\leq i,j\leq n+1},$
is given as follow,
\begin{equation}\label{eq omega n even}\small
\Omega=\left( {\begin{array}{*{20}c}
   0 & 0 &  \ldots   & 0 & {a}  \\
   0 & 0 &  \ldots   & {-a} & {\Omega _{2(n + 1)} }  \\
       \vdots  &  \vdots  &  \begin{array}{l}
 \,\,\,\,\,\,\,\,\,\,\,\,\,\,\qquad\qquad\qquad {\mathinner{\mkern2mu\raise1pt\hbox{.}\mkern2mu
 \raise4pt\hbox{.}\mkern2mu\raise7pt\hbox{.}\mkern1mu}}  \\
 \,\,\,\,\,\,\ \quad\quad \Omega_{ll}=(-1)^\frac{n}{2}a \\
  {\mathinner{\mkern2mu\raise1pt\hbox{.}\mkern2mu
 \raise4pt\hbox{.}\mkern2mu\raise7pt\hbox{.}\mkern1mu}}  \\
 \end{array}  &  \vdots  &  \vdots   \\
   0 & {-a} &  \ldots  & {\Omega_{nn}} & {\Omega _{n(n + 1)} }  \\
   {a} & {\Omega _{2(n + 1)} } &  \ldots  & {\Omega _{n(n + 1)} } & {\Omega_{(n+1)(n+1)}}  \\
\end{array}} \right),
\end{equation}
in which $l=\frac{n}{2}+1$.  Almost all
results of odd case are valid in even case. More precisely, we
can repeat Lemma \ref{lemm phi gm alpha}, Proposition \ref{prop
h} and Lemma \ref{lemm varthetaomega} exactly the same. But
for Lemma \ref{lemm y}, {\bf (i)} and {\bf (ii)} are valid, and {\bf (iii)} holds with some modification that we rewrite
it as follow.

\begin{lemm}\label{lemm y even}
The equation
$
\dot{S}=YS-\dot{z}S.\GM_\varpi,
$
 implies that
\begin{description}
  \item[(iii)] Moreover, if $S\Omega S^\tr=\Phi$, then
      $y_{i-1}=-y_{n-i}$, for $i=2,3,\ldots,\frac{n}{2}$.
      In the other word
$
Y\Phi=-\Phi Y^\tr.
$
\end{description}
\end{lemm}

Therefore to prove Theorem \ref{theo 1 int} in the even case, we can proceed analogously to the proof of the odd case.

\begin{rem}\rm
If we are more exact on the dimension of moduli space \Tm in
the even case and odd case, then we find a nice relationship
between them. Let \Tm be the moduli space given in Theorem \ref{theo 1 int} associated with a \cy $n$-fold,
where $n$ is even, and $\Tm'$ be the moduli space associated with a \cy
$(n-1)$-fold. Then we have
\[\small
\dim \Tm=\frac{n(n+2)}{4}+1=\frac{((n-1)+1)((n-1)+3)}{4}+1=\dim \Tm'.
\]
Thus, one of my interest for future works is to find more
relationships between structures of \Tm and $\Tm'$.
\end{rem}

\subsection{Five-Dimensional Case} \label{section 5-dim}
In this subsection we give an explicit presentation of \DHR vector
field $\H$, and in particular we verify its uniqueness by using
of self-duality of
Picard-Fuchs equation. Here we are following the notations and
terminologies of \S \ref{section cy n} for $n=5$. \\

The Picard-Fuchs equation (\ref{eq pf n 2-4})  associated with the fixed nowhere vanishing holomorphic $5$-form $\omega\in \F^5$ reduces to
\begin{equation} \label{eq pf 5}
\L=\vartheta^6-a_0(z)-a_1(z)\vartheta-a_2(z)\vartheta^2-a_3(z)\vartheta^3-a_4(z)\vartheta^4-a_5(z)\vartheta^5.
\end{equation}

\begin{lemm}\label{lemm coefficients}
The coefficients $a_i$'s of $\L$ given in {\rm (\ref{eq pf 5})} satisfy following equations:
\begin{equation}
\begin{array}{l}
a_3=-\frac{2}{3}a_4a_5+\frac{5}{3}a_5\vartheta a_5-\frac{5}{27}a_5^3-\frac{5}{3}\vartheta ^2a_5+2\vartheta a_4,\\ \\
a_1=\vartheta a_2-\vartheta^3a_4+\vartheta^4 a_5+\vartheta^2a_4a_5+\vartheta a_4\vartheta a_5+\frac{5}{3}a_5(\vartheta a_5)^2-\frac{1}{3}a_2a_5\\
\,\,\,\,\,\,\,\,+\frac{1}{27}a_4a_5^3
-\frac{10}{27}a_5^3\vartheta a_5+\frac{1}{81}a_5^5-\frac{1}{3}\vartheta a_4a_5^2-\frac{1}{3}a_4a_5\vartheta a_5+\frac{10}{9}a_5^2\vartheta^2a_5\\
\,\,\,\,\,\,\,\,
-\frac{10}{3}\vartheta a_5\vartheta^2a_5+\frac{1}{3}a_4\vartheta^2 a_5-\frac{5}{3}a_5\vartheta^3a_5.\\\\
\end{array}
\end{equation}
\end{lemm}
{\bf Proof. } By Proposition \ref{theo pp} the Picard-Fuchs equation (\ref{eq pf 5}) is self-dual, and the proof follows from Lemma \ref{lemm pp}{\bf (ii)} . \hfill\(\blacksquare\) \\

In the following proposition we compute all entries of the intersection matrix in the case $n=5$.

\begin{prop}\label{prop intersectionform 5}
The intersection form matrix $\Omega:=\left(\langle \omega_i,\omega_j \rangle\right)_{1\leq i,j\leq 6}$,
is given by
\begin{equation}\label{eq Omega 5}\small
\Omega=\left(
         \begin{array}{cccccc}
           0 & 0 & 0 & 0 & 0 & \tilde{a} \\
           0 & 0 & 0 & 0 & -\tilde{a} & \Omega_{26} \\
           0 & 0 & 0 & \tilde{a} & \Omega_{35} & \Omega_{36} \\
           0 & 0 & -\tilde{a} & 0 & \Omega_{45} & \Omega_{46} \\
           0 & \tilde{a} & -\Omega_{35} & -\Omega_{45} & 0 & \Omega_{56} \\
           -\tilde{a} & -\Omega_{26} & -\Omega_{36} & -\Omega_{46} & -\Omega_{56} & 0 \\
         \end{array}
       \right),
\end{equation}
where
$
\tilde{a}=c_0\exp \left( {\frac{1}{3} \int_0^z
a_5(v) \frac{dv}{v}}\right)
$
for some nonzero constant $c_0$, and
\[\small
\begin{array}{l}
\Omega_{26}=-\frac{2}{3}\tilde{a} a_5, \qquad\qquad\qquad\qquad\qquad
\Omega_{35}=\frac{1}{3}\tilde{a} a_5, \\\\
\Omega _{36} = \tilde a{a_4} + \frac{4}{9}\tilde aa_5^2 - \frac{2}{3}\tilde a\vartheta {a_5},\qquad\qquad\,\,
{\Omega _{45}} =  - \tilde a{a_4} - \frac{1}{3}\tilde aa_5^2 + \tilde a\vartheta {a_5},\\\\
{\Omega _{46}} =  - \tilde a{a_3} - \tilde a{a_4}{a_5} - \frac{8}{{27}}\tilde aa_5^3 + \tilde a\vartheta {a_4} + \frac{4}{3}\tilde a{a_5}\vartheta {a_5} - \frac{2}{3}\tilde a{\vartheta ^2}{a_5},\\\\
{\Omega _{56}} = \tilde a{a_2} + \frac{2}{3}\tilde a{a_3}{a_5} + \tilde aa_4^2 + \tilde a{a_4}a_5^2 + \frac{{16}}{{81}}\tilde aa_5^4 - \tilde a\vartheta {a_3} - \frac{5}{3}\tilde a{a_5}\vartheta {a_4}- \frac{{16}}{9}\tilde aa_5^2\vartheta {a_5}\\
\,\,\,\,\,\,\,\,\,\,\, - 2\tilde a{a_4}\vartheta {a_5} + \frac{4}{3}\tilde a{(\vartheta {a_5})^2} + \tilde a{\vartheta ^2}{a_4} + \frac{{16}}{9}\tilde a{a_5}{\vartheta ^2}{a_5}- \frac{2}{3}\tilde a{\vartheta ^3}{a_5}.
\end{array}\]
\end{prop}

{\bf Proof. } On account of (\ref{eq omega n}) we get that the matrix $\Omega$  is given by (\ref{eq Omega 5}), and we just need to find the entries $\Omega_{26},\Omega_{35},\Omega_{36},\Omega_{45},\Omega_{46},\Omega_{56}$. In order to do this, we first easily see that
$
\vartheta \tilde{a}=\frac{1}{3}\tilde{a}a_5.
$
By Picard-Fuchs equation (\ref{eq pf 5}) we obtain
\begin{equation}\label{eq pfo 5}
\vartheta\omega_6=\vartheta^6\omega=a_0\omega_1+a_1\omega_2+a_2\omega_3+a_3\omega_4+a_4\omega_5+a_5\omega_6.
\end{equation}
Since $\langle \omega_1,\omega_6 \rangle=\tilde{a}$, by considering (\ref{eq pfo 5}) and the fact that $\langle \omega_1,\omega_i \rangle=0,$ for $ i=1,2,\ldots,5,$ we find $\Omega_{26}$ as follow:
\begin{align}
\vartheta \langle \omega_1,\omega_6 \rangle&=\langle \omega_2,\omega_6 \rangle+\langle \omega_1,\vartheta\omega_6 \rangle \nonumber \\
&\Rightarrow \Omega_{26}=\vartheta\tilde{a}-\tilde{a}a_5 \Rightarrow \Omega_{26}=-\frac{2}{3}\tilde{a}a_5.\nonumber
\end{align}
We can find the rest of entries similarly and we just do that for $\Omega_{56}$,
\begin{align}
\Omega_{56}&=\vartheta \Omega_{46}-\langle \omega_4,\vartheta \omega_6 \rangle=\vartheta \Omega_{46}-a_2\Omega_{43}
           -a_4\Omega_{45}-a_5\Omega_{46} \nonumber\\
           &=\tilde a{a_2} + \frac{2}{3}\tilde a{a_3}{a_5} + \tilde aa_4^2 + \tilde a{a_4}a_5^2 + \frac{{16}}{{81}}\tilde aa_5^4
           - \tilde a\vartheta {a_3} - \frac{5}{3}\tilde a{a_5}\vartheta {a_4}- \frac{{16}}{9}\tilde aa_5^2\vartheta {a_5}\nonumber\\
           &- 2\tilde a{a_4}\vartheta {a_5} + \frac{4}{3}\tilde a{(\vartheta {a_5})^2} + \tilde a{\vartheta ^2}{a_4}
           + \frac{{16}}{9}\tilde a{a_5}{\vartheta ^2}{a_5}- \frac{2}{3}\tilde a{\vartheta ^3}{a_5}.  \nonumber
\end{align}   \hfill\(\blacksquare\)\\

The task is now to present \DHR vector field $\H$ explicitly. In order to do this, we use the chart $t$ that we pointed it out in Remark \ref{rem chart t}. In the theorem below, we verify that $\H$ as an ordinary differential equation is given as follow:
\begin{equation}\label{eq H 5}
\left \{ \begin{array}{l}
\dot t_0=\frac{t_0t_3}{t_1}
\\
\dot t_1=t_2
\\
\dot t_2=\frac{t_3^2t_4}{t_1t_6}
\\
\dot t_3=\frac{-t_2t_3t_6+t_3^2t_5}{t_1t_6}
\\
\dot t_4=\frac{\tilde{a}t_3t_6^2t_7}{t_1}
\\
\dot t_5=\frac{-t_3t_4+\tilde{a}t_3t_6^2t_8}{t_1}
\\
\dot t_6=\frac{-t_3t_5+\tilde{a}t_3t_6^2t_9}{t_1}
\\
\dot t_7=\frac{-t_3^2t_{10}}{t_1t_6}
\\
\dot t_8=\frac{-t_3^2t_{11}-t_3t_6t_7}{t_1t_6}
\\
\dot t_9=\frac{3\tilde{a}t_3t_5t_9-6\tilde{a}t_3t_6t_8-3t_3a_4-t_3a_5^2+3t_3\vartheta a_5}{3\tilde{a}t_1t_6}
\\
\dot t_{10}=-t_{12}
\\
\dot t_{11}=\frac{27\tilde{a}t_2t_{11}-54\tilde{a}t_3t_{10}+27\tilde{a}t_4t_8-27\tilde{a}t_5t_7+27a_2+27a_3a_5-27\vartheta a_3+15a_4a_5^2}{27\tilde{a}t_1}\\
\quad\,\,\,\, +\frac{-9a_4\vartheta a_5-54a_5\vartheta{a_4}+27\vartheta^2a_4+4a_5^4-42a_5^2\vartheta a_5+54a_5\vartheta^2a_5+18(\vartheta a_5)^2-18\vartheta^3a_5}{27\tilde{a}t_1}
\\
\dot t_{12}=\frac{-t_3a_0}{\tilde{a}t_1^2}
\end{array} \right. .
\end{equation}

\begin{theo} \label{teo 1 5}
Let \Tm be the moduli space introduced in Theorem \ref{theo 1 int},
for $n=5$. Then there is a chart $(t_0,t_1,\ldots,t_{12})$ for
$\Tm$ such that in this chart we obtain
 \[y_{1}=\frac{t_3^2}{t_1t_6},\qquad \& \qquad y_2=\frac{\tilde{a}t_3t_6^2}{t_1},\]
 and {\rm DHR} vector field $\H$ is given by (\ref{eq H 5}).
\end{theo}

{\bf Proof.} By Theorem \ref{theo 1 int} we get that \Tm is
13-dimensional. Using the equation
$
S\Omega S^\tr=\Phi,
$
and  considering $s_{11},s_{21},s_{22},s_{31},s_{32},s_{33},s_{41},s_{42},s_{43},s_{51},s_{52},s_{61}$ as independent entries of $S$, we thus can express  dependent entries in terms of independent entries as follows:
\begin{equation}\label{eq chart5}
\begin{array}{l}
s_{44}=\frac{1}{\tilde{a}s_{33}},\qquad\qquad\qquad\qquad\qquad
s_{53}=\frac{-3 \tilde{a}s_{32}s_{43}+3\tilde{a}s_{33}s_{42}+3a_4+a_5^2-3\vartheta a_5}{3\tilde{a}s_{22}},\\\\
s_{54}=\frac{-3s_{32}+s_{33}a_5}{3\tilde{a}s_{22}s_{33}},\qquad\qquad\qquad\quad
s_{55}=-\frac{1}{\tilde{a}s_{22}},\\\\
s_{62}=\frac{-27\tilde{a}s_{21}s_{52}+27\tilde{a}s_{22}s_{51}-27\tilde{a}s_{31}s_{42}+27\tilde{a}s_{32}s_{41}-27a_2-27a_3a_5+27\vartheta a_3-15a_4a_5^2}{27\tilde{a}s_{11}}\\
\quad\,\,\,\, +\frac{9a_4\vartheta a_5+54a_5\vartheta a_4-27\vartheta^2a_4-4a_5^4+42a_5^2\vartheta a_5-54a_5\vartheta^2 a_5-18(\vartheta a_5)^2+18 \vartheta^3 a_5}{27\tilde{a}s_{11}},\\\\
s_{63}=\frac{27\tilde{a}s_{21}s_{32}s_{43}-27\tilde{a}s_{21}s_{33}s_{42}-27\tilde{a}s_{22}s_{31}s_{43}+27\tilde{a}s_{22}s_{33}s_{41}-27s_{21}a_4-9s_{21}a_5^2}{27\tilde{a}s_{11}s_{22}}\\ \quad\,\,\,\, +\frac{27s_{21}\vartheta a_5-27s_{22}a_3-9s_{22}a_4a_5+27s_{22}\vartheta a_4-2s_{22}a_5^3+18s_{22}a_5\vartheta a_5-18s_{22}\vartheta^2a_5}{27\tilde{a}s_{11}s_{22}},\\\\
s_{64}=\frac{9s_{21}s_{32}-3s_{21}s_{33}a_5-9s_{22}s_{31}-9s_{22}s_{33}a_4-2s_{22}s_{33}a_5^2+6s_{22}s_{33}\vartheta a_5}{9\tilde{a}s_{11}s_{22}s_{33}},\\\\
s_{65}=\frac{3s_{21}-2s_{22}a_5}{3\tilde{a}s_{11}s_{22}},\qquad\qquad\qquad\quad
s_{66}=\frac{1}{\tilde{a}s_{11}}.
\end{array}
\end{equation}
We know that $\W$ is a family of one parameter 5-dimensional
Calabi-Yau manifolds parameterized by $z$. Hence we present the chart $t=(t_0,t_1,\ldots,t_{12})$ for $\Tm$, where
$t_0=z,\,\,\, t_1=s_{11},\,\,\,
t_2=s_{21},\,\,\,t_3=s_{22},\,\,\, t_4=s_{31}, \,\,\,
t_5=s_{32},\,\,\, t_6=s_{33},\,\,\, t_7=s_{41},\,\,\,
t_8=s_{42},\,\,\, t_9=s_{43},\,\,\, t_{10}=s_{51},\,\,\,
t_{11}=s_{52},\,\,\, t_{12}=s_{61}$.
The same as the proof of Proposition \ref{prop h},
 define $\H:=\sum_{i=0}^{12}\H_i(t)\frac{\partial}{\partial t_i}$
and  set $\dot{t}_i:=\H_i(t)$. Then $\H_0$, $y_1$ and $y_2$ follow from Lemma \ref{lemm y}. Therefore, the right hand side of the equation
\begin{equation}\label{eq syh5}
\dot{S}=YS-\dot{z}S.\widehat{\GM}_\varpi,
\end{equation}
is totally determined. We also substitute $a_1$, $a_3$ from Lemma \ref{lemm
coefficients}, and dependent entries from (\ref{eq chart5}) in the right hand side of (\ref{eq syh5}).
Consequently, the rest of $\H_i$'s follow directly from (\ref{eq syh5}). Note that
the compatibility of $S\Omega S^\tr=\Phi$ and $\dot{S}=YS-\dot{z}S.\widehat{\GM}_\varpi$ follow from
substituting $a_1$ and $a_3$ from Lemma \ref{lemm coefficients}.
\hfill\(\blacksquare\)

\subsection{Three-Dimensional Case}\label{section 3-dim}
Here we substitute the dimension $n=5$ with $n=3$, and  will proceed analogously to  \S\ref{section 5-dim}.

The Picard-Fuchs equation $\L$ is given as
\begin{equation} \label{eq pf 3}
\L=\vartheta^4-a_0(z)-a_1(z)\vartheta-a_2(z)\vartheta^2-a_3(z)\vartheta^3,
\end{equation}
where coefficients $a_i$'s, by Lemma \ref{lemm pp}{\bf (i)}, satisfy the following relationship:
\begin{equation}
a_1=\frac{3}{4}a_3\vartheta a_3+\vartheta a_2-\frac{1}{2}\vartheta^2a_3-\frac{1}{8}a_3^3-\frac{1}{2}a_2a_3.
\end{equation}
The intersection form matrix is as follow
\begin{equation}\label{eq Omega 3}
\Omega=\left(
         \begin{array}{cccc}
           0 & 0 & 0 & \tilde{a} \\
           0 & 0 & -\tilde{a} & \Omega_{24} \\
           0 & \tilde{a} & 0 & \Omega_{34} \\
           -\tilde{a} & -\Omega_{24} & -\Omega_{34} & 0 \\
         \end{array}
       \right),
\end{equation}
in which
\begin{equation}\label{eq atild 3}
\tilde{a}=c_0\exp \left( {\frac{1}{2} \int_0^za_3(v) \frac{dv}{v}}\right),
\end{equation}
and
\[\begin{array}{l}
\Omega_{24}=-\frac{1}{2}\tilde{a} a_3,\qquad \&\qquad
\Omega_{34}=\frac{1}{4}\tilde{a} a_3^2+\tilde{a}a_2-\frac{1}{2}\tilde{a}\vartheta a_3.
\end{array}\]
The chart $t=(t_0,t_1,\ldots,t_{6})$ of
$\Tm$ is obtained by settin $t_0=z,\,\,\,
t_1=s_{11},\,\,\, t_2=s_{21},\,\,\,t_3=s_{22},\,\,\,
t_4=s_{31}, \,\,\, t_5=s_{32},\,\,\, t_6=s_{41}$. We compute below dependent entries:
\begin{equation}\label{eq chart3}
\begin{array}{l}
s_{33}=-\frac{1}{\tilde{a}s_{22}},\qquad\qquad\qquad\qquad
s_{42}=\frac{4 \tilde{a}s_{22}s_{31}-4\tilde{a}s_{21}s_{32}-a_3^2-4a_2+2\vartheta a_3}{4\tilde{a}s_{11}},\\\\
s_{43}=\frac{2s_{21}-s_{22}a_3}{2\tilde{a}s_{11}s_{22}},\qquad\qquad\qquad\,\,
s_{44}=\frac{1}{\tilde{a}s_{11}}.
\end{array}
\end{equation}
In the chart $t$, the meromorphic function $y_1$ is given by $$y_{1}=-\frac{\tilde{a}t_3^3}{t_1},$$
 and \DHR vector field $\H$ has following presentation:
\begin{equation}
\label{eq H 3}
\left \{ \begin{array}{l}
\dot t_0=\frac{t_0t_3}{t_1}
\\
\dot t_1=t_2

\\
\dot t_2=-\frac{\tilde{a}t_3^3t_4}{t_1}
\\
\dot t_3=-\frac{t_2t_3+\tilde{a}t_3^3t_5}{t_1}

\\
\dot t_4=-t_6
\\
\dot t_5=\frac{4\tilde{a}t_2t_5-8\tilde{a}t_3t_4+a_3^2+4a_2-2\vartheta a_3}{4\tilde{a}t_1}
\\
\dot t_6=-\frac{t_3a_0}{\tilde{a}t_1^2}

\end{array} \right. .
\end{equation}

\begin{exam}
In Example \ref{exam 14 cases} we introduced 14 families of \cy 3-folds for which Theorem \ref{theo 1 int} holds. We can rewrite the Picard-Fuchs equation (\ref{eq pf3 14}) as follow
\[
\L=\vartheta^4-a_0(z)-a_1(z)\vartheta-a_2(z)\vartheta^2-a_3(z)\vartheta^3,
\]
in which
\begin{align}
&a_0(z)=\frac{cz(r_1r_2-r_1^2r_2-r_1r_2^2+r_1^2r_2^2)}{1-cz}, \qquad  \qquad a_1(z)=\frac{cz(r_1+r_2-r_1^2-r_2^2)}{1-cz}, \nonumber\\
&a_2(z)=\frac{cz(1+r_1+r_2-r_1^2-r_2^2)}{1-cz}, \qquad \qquad \qquad a_3(z)=\frac{2cz}{1-cz}. \nonumber
\end{align}
On account of (\ref{eq atild 3}), we obtain
\[
\tilde{a}=\frac{c_0}{1-cz}.
\]
Now by replacing $a_i$'s and $\tilde{a}$ in (\ref{eq H 3}) we find \DHR vector field $\H$ on the moduli space associated with any of the families given in Table \ref{table 14 case}.
\end{exam}




\end{document}